\documentclass[12pt,twoside,british]{article}
\usepackage[T1]{fontenc}
\usepackage[utf8]{inputenc}
\usepackage[a4paper]{geometry}
\usepackage{color}
\usepackage{float}
\usepackage{amsmath}
\usepackage{amsthm}
\usepackage{amssymb}
\usepackage{graphicx}
\usepackage{setspace}

\makeatletter

\providecommand{\tabularnewline}{\\}

\numberwithin{equation}{section}
\numberwithin{figure}{section}
\theoremstyle{plain}
\newtheorem{thm}{\protect\theoremname}[section]
\theoremstyle{definition}
\newtheorem{defn}[thm]{\protect\definitionname}
\newcommand{\lyxaddress}[1]{
	\par {\raggedright #1
	\vspace{1.4em}
	\noindent\par}
}

\usepackage{amsthm}
\usepackage{setspace}

\numberwithin{equation}{section}
\numberwithin{figure}{section}

\usepackage{amsthm}
\usepackage{eucal}
\usepackage{epstopdf}
\usepackage{graphicx}
\theoremstyle{plain}
\DeclareMathOperator*{\argmin}{arg\,min} %argmin
\newtheoremstyle{boldremark}
    {\dimexpr\topsep/2\relax} % space above
    {\dimexpr\topsep/2\relax} % space below
    {}          % body font
    {}          % indent amount
    {\bfseries} % theorem head font
    {.}         % punctuation after theorem head
    {.5em}      % space after theorem head
    {}          % theorem hed spec. (empty = "normal")

\theoremstyle{boldremark}

\usepackage{fancyhdr}
\usepackage{blindtext}
\usepackage{titlesec}

\titleformat{name=\chapter}[display]
{\normalfont\Huge\itshape}
{\titlerule[1pt]\vspace{-33pt}\filleft%
  \parbox[t]{6em}{%
    \raggedleft%
    \rule{\linewidth}{0.5ex}\newline%
    \chaptertitlename\ \thechapter%
  }%
}
{4pc}
{\normalfont\upshape\bfseries\Huge}
\allowdisplaybreaks

\usepackage{babel}

\makeatother

\usepackage{babel}
\addto\captionsbritish{\renewcommand{\definitionname}{Definition}}
\addto\captionsbritish{\renewcommand{\theoremname}{Theorem}}
\addto\captionsenglish{\renewcommand{\definitionname}{Definition}}
\addto\captionsenglish{\renewcommand{\theoremname}{Theorem}}
\addto\captionsitalian{\renewcommand{\definitionname}{Definizione}}
\addto\captionsitalian{\renewcommand{\theoremname}{Teorema}}
\providecommand{\definitionname}{Definition}
\providecommand{\theoremname}{Theorem}

\usepackage[colorlinks,pdfpagelabels,pdfstartview = FitH,bookmarksopen
= true,bookmarksnumbered = true,linkcolor = red,plainpages =
false,hypertexnames = false,citecolor = blue,pagebackref=false]{hyperref}

\begin{document}
\title{\textbf{A physics-informed deep learning approach for\\solving strongly
degenerate parabolic problems}}
\author{Pasquale Ambrosio, Salvatore Cuomo, Mariapia De Rosa}
\date{February 15, 2024}
\maketitle
\begin{abstract}
\noindent In recent years, Scientific Machine Learning (SciML) methods
for solving Partial Differential Equations (PDEs) have gained increasing
popularity. Within such a paradigm, Physics-Informed Neural Networks
(PINNs) are novel deep learning frameworks for solving initial-boundary
value problems involving nonlinear PDEs. Recently, PINNs have shown
promising results in several application fields. Motivated by applications
to gas filtration problems, here we present and evaluate a PINN-based
approach to predict solutions to \textit{strongly degenerate parabolic
problems with asymptotic structure of Laplacian type}. To the best
of our knowledge, this is one of the first papers demonstrating the
efficacy of the PINN framework for solving such kind of problems.
In particular, we estimate an appropriate approximation error for
some test problems whose analytical solutions are fortunately known.
The numerical experiments discussed include two and three-dimensional
spatial domains, emphasizing the effectiveness of this approach in
predicting accurate solutions.\vspace{0.15cm}
\end{abstract}

\section{Introduction}

\noindent $\hspace*{1em}$In this paper, we aim to exploit a novel
Artificial Intelligence (AI) methodology, known as \textit{Physics-Informed
Neural Networks} (PINNs), to predict solutions to Cauchy-Dirichlet
problems of the type 
\begin{equation}
\begin{cases}
\begin{array}{cc}
\partial_{t}u-\mathrm{div}\left((\vert\nabla u\vert-1)_{+}\,\,\frac{\nabla u}{\vert\nabla u\vert}\right)=f & \,\,\,\mathrm{in}\,\,\,\Omega_{T}:=\Omega\times(0,T),\\
u=w & \mathrm{on}\,\,\,\partial_{\mathrm{par}}\Omega_{T},\,\,\,\,\,\,\,\,\,\,\,\,\,\,\,\,\,\,\,\,\,\,
\end{array}\end{cases}\label{eq:mainprob}
\end{equation}
where $\Omega$ is a bounded connected open subset of $\mathbb{R}^{n}$
($2\leq n\leq3$) with Lipschitz boundary, $f$ and $w$ are given
real-valued functions defined over $\overline{\Omega}\times[0,T]$
and $\partial_{\mathrm{par}}\Omega_{T}$ respectively, $\nabla u$
denotes the spatial gradient of an unknown solution $u:\overline{\Omega}\times[0,T)\rightarrow\mathbb{R}$,
while $\left(\,\cdot\,\right)_{+}$ stands for the positive part.\\
$\hspace*{1em}$A motivation for studying problem (\ref{eq:mainprob})
can be found in \textit{gas filtration problems} (see \cite{Akh}
and \cite{AmbPass}). In order to make the paper self-contained, we
provide a brief explanation in Section \ref{subsec:Motivation} below.\\
$\hspace*{1em}$As for the parabolic equation $(\ref{eq:mainprob})_{1}$,
the regularity properties of its weak solutions have been recently
studied in \cite{Amb,AmbPass} and \cite{GenPass}. The main novelty
of this PDE is that it exhibits a strong degeneracy, coming from the
fact that its modulus of ellipticity vanishes in the region $\{\vert\nabla u\vert\leq1\}$,
and hence its principal part behaves like a \textit{Laplace operator
only at infinity}.\\
$\hspace*{1em}$The regularity of solutions to parabolic problems
with asymptotic structure of Laplacian type had already been investigated
in \cite{Ise}, where a BMO\footnote{BMO denotes the space of functions with \textit{bounded mean oscillations} (see \cite[Chapter 2]{Giu}).} regularity was proved for solutions to
asymptotically parabolic systems in the case $f=0$ (see also \cite{KuuMing},
where the local Lipschitz continuity of weak solutions with respect
to the spatial variable is established). In addition, we want to mention
the results contained in \cite{Byun}, where nonhomogeneous parabolic
problems with an asymptotic regularity in divergence form of $p$-Laplacian
type are considered. There, Byun, Oh and Wang establish a global Calderón-Zygmund
estimate by converting a given asymptotically regular problem to a
suitable regular problem.\\
$\hspace*{1em}$Concerning the approach used here, the PINNs are a
Scientific Machine Learning (SciML) technique based on Artificial Neural
Networks (ANNs) with the feature of adding constraints to make the
predicted results more in line with the physical laws of the addressed
problem. The concept of PINNs was introduced
in \cite{Karn,Raissi1,Raissi2} and \cite{Raissi3}
to solve PDE-based problems.
The PINNs predict the solution to a PDE under prescribed initial-boundary
conditions by training a neural network to minimize a cost function,
called \textit{loss function}, which penalizes some suitable terms
on a set of admissible functions $u$ (for more information, we refer the interested reader to \cite{Cuomo2}).\\
$\hspace*{1em}$The kind of approach we want to propose here can offer
effective solutions to real problems such as (\ref{eq:mainprob})
and can be applied in many other different fields: for example, in
production and advanced engineering \cite{Zobeiry}, for transportation
problems \cite{Fraces}, and for virtual thermal sensors using real-time
simulations \cite{Go}. Additionally, it is employed to solve groundwater
flow equations \cite{Cuomo} and address petroleum and gas contamination
\cite{Soriano}.\\
$\hspace*{1em}$As far as we know, this is one of the first papers
demonstrating the effectiveness of the PINN framework for solving
strongly degenerate parabolic problems of the type (\ref{eq:mainprob}).\smallskip{}

\subsection{Motivation}\label{subsec:Motivation}
\noindent $\hspace*{1em}$Before describing the structure of this
paper, we wish to motivate our study by pointing out that, in the
physical cases $n=2$ and $n=3$, degenerate equations of the form
$(\ref{eq:mainprob})_{1}$ may arise in \textit{gas filtration problems
taking into account the initial pressure gradient}. \\
$\hspace*{1em}$The existence of remarkable deviations from the linear
Darcy filtration law has been observed in several systems consisting
of a fluid and a porous medium (e.g., the filtration of a gas in argillous
rocks). One of the manifestations of this nonlinearity is the existence
of a limiting pressure gradient, i.e. the minimum value of the pressure
gradient for which fluid motion takes place. In general, fluid motion
still occurs for subcritical values of the pressure gradient, but
very slowly; when achieving the limiting value of the pressure gradient,
there is a marked acceleration of the filtration. Therefore, the limiting-gradient
concept provides a good approximation for velocities that are not
too low.\\
$\hspace*{1em}$In accordance with some experimental results (see
\cite{Akh}), under certain physical conditions one can take the gas
filtration law in the very simple form 
\[
\begin{cases}
\begin{array}{cc}
\mathbf{v}=-\,\frac{k}{\mu}\,\nabla p\left[1-\frac{\beta}{\left|\nabla p^{2}\right|}\right] & \,\,\,\mathrm{if}\,\,\left|\nabla p^{2}\right|\geq\beta,\\
\mathbf{v}=\mathbf{0}\,\,\,\,\,\,\,\,\,\,\,\,\,\,\,\,\,\,\,\,\,\,\,\,\,\,\,\,\,\,\,\,\,\,\,\,\,\,\,\,\,\,\, & \,\,\,\mathrm{if}\,\,\left|\nabla p^{2}\right|<\beta,
\end{array}\end{cases}
\]
where $\mathbf{v}=\mathbf{v}(x,t)$ is the filtration velocity, $k$
is the rock permeability, $\mu$ is the gas viscosity, $p=p(x,t)$
is the pressure and $\beta$ is a positive constant. Under this assumption
we obtain a particularly simple expression for the gas mass velocity
(flux) $\boldsymbol{j}$, which contains only the gradient of the
pressure squared, exactly as in the usual gas filtration problems:
\begin{equation}
\begin{cases}
\begin{array}{cc}
\boldsymbol{j}=\varrho\,\mathbf{v}=-\,\frac{k}{2\mu C}\left[\nabla p^{2}-\beta\,\frac{\nabla p^{2}}{\left|\nabla p^{2}\right|}\right] & \,\,\,\mathrm{if}\,\,\left|\nabla p^{2}\right|>\beta,\\
\boldsymbol{j}=\mathbf{0}\,\,\,\,\,\,\,\,\,\,\,\,\,\,\,\,\,\,\,\,\,\,\,\,\,\,\,\,\,\,\,\,\,\,\,\,\,\,\,\,\,\,\,\,\,\,\,\,\,\,\,\,\,\,\,\,\,\,\,\,\,\,\,\,\,\,\,\,\,\, & \,\,\,\mathrm{if}\,\,\left|\nabla p^{2}\right|\leq\beta,
\end{array}\end{cases}\label{eq:flux}
\end{equation}
where $\varrho$ is the gas density and $C$ is a positive constant.
Plugging expression (\ref{eq:flux}) into the gas mass-conservation
equation, we obtain the basic equation for the pressure: 
\begin{equation}
\begin{cases}
\begin{array}{cc}
\frac{\partial p}{\partial t}=\,\frac{k}{2m\mu}\,\mathrm{div}\left[\nabla p^{2}-\beta\,\frac{\nabla p^{2}}{\left|\nabla p^{2}\right|}\right] & \,\,\,\mathrm{if}\,\,\left|\nabla p^{2}\right|>\beta,\\
\,\,\frac{\partial p}{\partial t}=0\,\,\,\,\,\,\,\,\,\,\,\,\,\,\,\,\,\,\,\,\,\,\,\,\,\,\,\,\,\,\,\,\,\,\,\,\,\,\,\,\,\,\,\,\,\,\,\,\,\,\,\,\,\,\,\,\,\,\,\, & \,\,\,\mathrm{if}\,\,\left|\nabla p^{2}\right|\leq\beta,
\end{array}\end{cases}\label{eq:pressure}
\end{equation}
where $m$ is a positive constant. Equation (\ref{eq:pressure}) implies,
first of all, that the steady gas motion is described by the same
relations as in the steady motion of an incompressible fluid if we
replace the pressure of the incompressible fluid with the square of
the gas pressure. In addition, if the gas pressure differs very little
from some constant pressure $p_{0}$, or if the gas pressure differs
considerably from a constant value only in regions where the gas motion
is nearly steady, then the equation for the gas filtration in the
region of motion can be ``linearized'' following L. S. Leibenson,
and thus obtaining (see \cite{Akh} again)
\begin{equation}
\begin{cases}
\begin{array}{cc}
\frac{\partial p^{2}}{\partial t}=\,\frac{k\,p_{0}}{m\mu}\,\mathrm{div}\left[\nabla p^{2}-\beta\,\frac{\nabla p^{2}}{\left|\nabla p^{2}\right|}\right] & \,\,\,\mathrm{if}\,\,\left|\nabla p^{2}\right|>\beta,\\
\,\,\frac{\partial p^{2}}{\partial t}=0\,\,\,\,\,\,\,\,\,\,\,\,\,\,\,\,\,\,\,\,\,\,\,\,\,\,\,\,\,\,\,\,\,\,\,\,\,\,\,\,\,\,\,\,\,\,\,\,\,\,\,\,\,\,\,\,\,\,\,\, & \,\,\,\mathrm{if}\,\,\left|\nabla p^{2}\right|\leq\beta.
\end{array}\end{cases}\label{eq:sqpress}
\end{equation}
$\hspace*{1em}$Setting $u=p^{2}$ and performing a suitable scaling,
equation (\ref{eq:sqpress}) turns into 
\[
\frac{\partial u}{\partial t}-\mathrm{div}\left((\vert\nabla u\vert-1)_{+}\,\,\frac{\nabla u}{\vert\nabla u\vert}\right)=0,
\]
which is nothing but equation $(\ref{eq:mainprob})_{1}$ with $f\equiv0$.
This is why $(\ref{eq:mainprob})_{1}$ is sometimes called the \textit{Leibenson
equation} in the literature.\\
\\
$\hspace*{1em}$The paper is organized as follows. Section \ref{sec:preliminaries}
is devoted to the preliminaries: after a list of some classic notations,
we provide details on the strongly degenerate parabolic problem (\ref{eq:mainprob}).
In Section \ref{sec:method}, we describe the PINN methodology that was employed. Section \ref{sec:results}
presents the results that were obtained. Finally, Section \ref{sec:conclusions}
provides the conclusions.\vspace{0.6cm}

\section{Notation and preliminaries}\label{sec:preliminaries}
\noindent $\hspace*{1em}$In what follows, the
norm we use on $\mathbb{R}^{n}$ will be the standard Euclidean one
and it will be denoted by $\left|\,\cdot\,\right|$. In particular,
for the vectors $\xi,\eta\in\mathbb{R}^{n}$, we write $\langle\xi,\eta\rangle$
for the usual inner product and $\left|\xi\right|:=\langle\xi,\xi\rangle^{\frac{1}{2}}$
for the corresponding Euclidean norm. For points in space-time, we
will frequently use abbreviations like $z=(x,t)$ or $z_{0}=(x_{0},t_{0})$,
for spatial variables $x$, $x_{0}\in\mathbb{R}^{n}$ and times $t$,
$t_{0}\in\mathbb{R}$. We also denote by $B_{\rho}(x_{0})=\left\{ x\in\mathbb{R}^{n}:\left|x-x_{0}\right|<\rho\right\} $
the open ball with radius $\rho>0$ and center $x_{0}\in\mathbb{R}^{n}$.
Moreover, we use the notation
\[
Q_{\rho}(z_{0}):=B_{\rho}(x_{0})\times(t_{0}-\rho^{2},t_{0}),\,\,\,\,\,z_{0}=(x_{0},t_{0})\in\mathbb{R}^{n}\times\mathbb{R},\,\,\rho>0,
\]
for the backward parabolic cylinder with vertex $(x_{0},t_{0})$ and
width $\rho$. Finally, for a general cylinder $Q=A\times(t_{1},t_{2})$,
where $A\subset\mathbb{R}^{n}$ and $t_{1}<t_{2}$, we denote by
\[
\partial_{\mathrm{par}}Q:=(\overline{A}\times\{t_{1}\})\cup(\partial A\times(t_{1},t_{2}))
\]
the usual \textit{parabolic boundary} of $Q$.\\
$\hspace*{1em}$To give the definition of a weak solution to problem
$(\ref{eq:mainprob})$, we now introduce the function $H:\mathbb{R}^{n}\rightarrow\mathbb{R}^{n}$
defined by 
\[
H(\xi):=\begin{cases}
\begin{array}{cc}
(\left|\xi\right|-1)_{+}\,\frac{\xi}{\left|\xi\right|} & \,\,\,\,\,\,\,\,\mathrm{if}\,\,\,\xi\in\mathbb{R}^{n}\,\backslash\,\{0\},\\
0 & \mathrm{if}\,\,\,\xi=0.\,\,\,\,\,\,\,\,\,\,
\end{array}\end{cases}
\]

\begin{defn}
\noindent Let $f\in L_{loc}^{1}(\Omega_{T})$. A function $u\in C^{0}\left((0,T);L^{2}(\Omega)\right)\cap L^{2}\left(0,T;W^{1,2}(\Omega)\right)$
is a \textit{weak solution} of equation $(\ref{eq:mainprob})_{1}$
if and only if for any test function $\varphi\in C_{0}^{\infty}(\Omega_{T})$
the following integral identity holds:
\begin{equation}
\int_{\Omega_{T}}\left(u\cdot\partial_{t}\,\varphi-\langle H(\nabla u),\nabla\varphi\rangle\right)\,dz\,=\,-\int_{\Omega_{T}}f\varphi\,dz.\label{eq:weak}
\end{equation}
\end{defn}

\noindent 
\begin{defn}
\noindent Let $w\in C^{0}\left([0,T];L^{2}(\Omega)\right)\cap L^{2}\left(0,T;W^{1,2}(\Omega)\right)$.
We identify a function
\[
u\in C^{0}\left([0,T];L^{2}(\Omega)\right)\cap L^{2}\left(0,T;W^{1,2}(\Omega)\right)
\]
as a \textit{weak solution of the Cauchy-Dirichlet problem }(\ref{eq:mainprob})
if and only if (\ref{eq:weak}) holds and, moreover, $u\in w+L^{2}\left(0,T;W_{0}^{1,2}(\Omega)\right)$
and $u(\cdot,0)=w(\cdot,0)$ in the $L^{2}$-sense, that is 
\begin{equation}
\underset{t\,\searrow\,0}{\lim}\,\Vert u(\cdot,t)-w(\cdot,0)\Vert_{L^{2}(\Omega)}\,=\,0.\label{eq:L2sense}
\end{equation}
Therefore, the initial condition $u=w$ on $\overline{\Omega}\times\{0\}$
has to be understood in the usual $L^{2}$-sense (\ref{eq:L2sense}),
while the condition $u=w$ on the lateral boundary $\partial\Omega\times(0,T)$
has to be meant in the sense of traces, i.e. $(u-w)\,(\cdot,t)\in W_{0}^{1,2}(\Omega)$
for almost every $t\in(0,T)$.
\end{defn}

\noindent $\hspace*{1em}$Taking $p=2$ and $\nu=1$ in \cite[Theorem 1.1]{AmbPass},
we immediately obtain the following spatial Sobolev regularity result:
\begin{thm}
\noindent \label{thm:theo1} Let $n\geq2$, $\frac{2n\,+\,4}{n\,+\,4}\leq q<\infty$
and $f\in L^{q}\left(0,T;W^{1,q}(\Omega)\right)$. Moreover, assume
that 
\[
u\in C^{0}\left((0,T);L^{2}(\Omega)\right)\cap L^{2}\left(0,T;W^{1,2}(\Omega)\right)
\]
is a weak solution of equation $\mathrm{(\ref{eq:mainprob})_{1}}$.
Then the solution satisfies 
\[
H(\nabla u)\,\in\,L_{loc}^{2}\left(0,T;W_{loc}^{1,2}(\Omega,\mathbb{R}^{n})\right).
\]
Furthermore, the following estimate
\[
\int_{Q_{\rho/2}(z_{0})}\left|\nabla H(\nabla u)\right|^{2}dz\,\leq\,\,c\left(\Vert\nabla f\Vert_{L^{q}(Q_{R_{0}})}+\,\Vert\nabla f\Vert_{L^{q}(Q_{R_{0}})}^{2}\right)\,+\,\frac{c}{R^{2}}\left(\Vert\nabla u\Vert_{L^{2}(Q_{R_{0}})}^{2}+\,1\right)
\]
holds true for any parabolic cylinder $Q_{\rho}(z_{0})\subset Q_{R}(z_{0})\subset Q_{R_{0}}(z_{0})\Subset\Omega_{T}$
and a positive constant $c$ depending on $n$, $q$ and $R_{0}$.
\end{thm}

\noindent $\hspace*{1em}$From the above result one can easily deduce
that $u$ admits a weak time derivative $u_{t}$, which belongs to
the local Lebesgue space $L_{loc}^{\min\,\{2,\,q\}}(\Omega_{T})$.
The idea is roughly as follows. Consider equation $(\ref{eq:mainprob})_{1}$;
since the previous theorem tells us that in a certain pointwise sense
the second spatial derivatives of $u$ exist, we may develop the expression
under the divergence symbol; this will give us an expression that
equals $u_{t}$, from which we get the desired summability of the
time derivative. Such an argument has been made rigorous in \cite[Theorem 1.2]{AmbPass},
from which we can derive the next result.

\begin{thm}
\noindent \label{thm:theo2} Under the assumptions of Theorem \ref{thm:theo1},
the time derivative of the solution exists in the weak sense and satisfies
\[
\partial_{t}u\,\in\,L_{loc}^{\min\,\{2,\,q\}}(\Omega_{T}).
\]
Furthermore, the following estimate\begin{align*}
\left(\int_{Q_{\rho/2}(z_{0})}\left|\partial_{t}u\right|^{\min\,\{2,\,q\}}\,dz\right)^{\frac{1}{\min\,\{2,\,q\}}}&\leq\,\,c\,\Vert f\Vert_{L^{q}(Q_{R_{0}})}+\,c\,\left(\Vert \nabla f\Vert_{L^{q}(Q_{R_{0}})}+\,\Vert \nabla f\Vert_{L^{q}(Q_{R_{0}})}^{2}\right)^{\frac{1}{2}}\\
&\,\,\,\,\,\,\,+\frac{c}{R}\left(\Vert \nabla u\Vert_{L^{2}(Q_{R_{0}})}^{2}+1\right)^{\frac{1}{2}}
\end{align*}

\noindent holds true for any parabolic cylinder $Q_{\rho}(z_{0})\subset Q_{R}(z_{0})\subset Q_{R_{0}}(z_{0})\Subset\Omega_{T}$
and a positive constant $c$ depending on $n$, $q$ and $R_{0}$. 
\end{thm}

\noindent $\hspace*{1em}$Now, let the assumptions of Theorem \ref{thm:theo1}
be in force. For $\varepsilon\in[0,1]$ and a couple of standard,
non-negative, radially symmetric mollifiers $\phi_{1}\in C_{0}^{\infty}(B_{1}(0))$
and $\phi_{2}\in C_{0}^{\infty}((-1,1))$ we define 
\[
f_{\varepsilon}(x,t):=\int_{-1}^{1}\int_{B_{1}(0)}f(x-\varepsilon y,t-\varepsilon s)\,\phi_{1}(y)\,\phi_{2}(s)\,dy\,ds,
\]
where $f$ is meant to be extended by zero outside $\Omega_{T}$.
Observe that $f_{0}=f$ and $f_{\varepsilon}\in C^{\infty}(\Omega_{T})$
for every $\varepsilon\in(0,1]$.\\
Next, we consider a domain in space-time denoted by $\Omega'_{1,2}:=\Omega'\times(t_{1},t_{2})$,
where $\Omega'\subseteq\Omega$ is a bounded domain with smooth boundary
and $(t_{1},t_{2})\subseteq(0,T)$. In the following, we will need
the definitions below.
\begin{defn}
\noindent Let $\varepsilon\in(0,1]$. A function $u_{\varepsilon}\in C^{0}\left((t_{1},t_{2});L^{2}(\Omega')\right)\cap L^{2}\left(t_{1},t_{2};W^{1,2}(\Omega')\right)$
is a \textit{weak solution} of the equation
\begin{equation}
\partial_{t}\,u_{\varepsilon}-\mathrm{div}\left(H(\nabla u_{\varepsilon})+\varepsilon\,\nabla u_{\varepsilon}\right)=f_{\varepsilon}\,\,\,\,\,\,\,\,\mathrm{in}\,\,\,\Omega'_{1,2}\label{eq:regequ}
\end{equation}
if and only if for any test function $\varphi\in C_{0}^{\infty}(\Omega'_{1,2})$
the following integral identity holds:
\begin{equation}
\int_{\Omega'_{1,2}}\left(u_{\varepsilon}\cdot\partial_{t}\,\varphi-\langle H(\nabla u_{\varepsilon})+\varepsilon\,\nabla u_{\varepsilon},\nabla\varphi\rangle\right)\,dz\,=\,-\int_{\Omega'_{1,2}}f_{\varepsilon}\,\varphi\,dz.\label{eq:weak2}
\end{equation}
\end{defn}

\noindent 
\begin{defn}
\noindent \label{def:L2-traces} Let $\varepsilon\in(0,1]$ and $u\in C^{0}\left([t_{1},t_{2}];L^{2}(\Omega')\right)\cap L^{2}\left(t_{1},t_{2};W^{1,2}(\Omega')\right)$.
We identify a function 
\[
u_{\varepsilon}\in C^{0}\left([t_{1},t_{2}];L^{2}(\Omega')\right)\cap L^{2}\left(t_{1},t_{2};W^{1,2}(\Omega')\right)
\]
as a \textit{weak solution of the Cauchy-Dirichlet problem\medskip{}
}
\begin{equation}
\begin{cases}
\begin{array}{cc}
\partial_{t}\,u_{\varepsilon}-\mathrm{div}\left(H(\nabla u_{\varepsilon})+\varepsilon\,\nabla u_{\varepsilon}\right)=f_{\varepsilon} & \mathrm{in}\,\,\Omega'_{1,2},\,\,\,\,\,\,\,\,\,\,\,\\
u_{\varepsilon}=u & \,\,\,\,\,\,\mathrm{on}\,\,\partial_{\mathrm{par}}\Omega'_{1,2},\,\,\,\,\,\,\,
\end{array}\end{cases}\label{eq:CAUCHYDIR}
\end{equation}
if and only if (\ref{eq:weak2}) holds and, moreover, 
\[
u_{\varepsilon}\,\in\,u+L^{2}\left(t_{1},t_{2};W_{0}^{1,2}(\Omega')\right),
\]
$u_{\varepsilon}(\cdot,t_{1})=u(\cdot,t_{1})$ in the usual $L^{2}$-sense
and the condition $u_{\varepsilon}=u$ on the lateral boundary $\partial\Omega'\times(t_{1},t_{2})$
holds in the sense of traces, i.e. $(u_{\varepsilon}-u)\,(\cdot,t)\in W_{0}^{1,2}(\Omega')$
for almost every $t\in(t_{1},t_{2})$. 
\end{defn}

\noindent $\hspace*{1em}$Due to the strong degeneracy of equation
$(\ref{eq:mainprob})_{1}$, in order to prove Theorems \ref{thm:theo1}
and \ref{thm:theo2} above, the authors of \cite{AmbPass} resort
to the family of approximating parabolic problems (\ref{eq:CAUCHYDIR}).
These problems exhibit a milder degeneracy than $(\ref{eq:mainprob})$
and the advantage of considering them stems from the fact that the
existence of a unique energy solution $u_{\varepsilon}$ satisfying
the requirements of Definition \ref{def:L2-traces} can be ensured
by the classic existence theory for parabolic equations (see \cite[Chapter 2, Theorem 1.2 and Remark 1.2]{Lions}).\\
$\hspace*{1em}$Moreover, if $\Omega'_{1,2}=Q_{R_{0}}(z_{0}):=B_{R_{0}}(x_{0})\times(t_{0}-R_{0}^{2},t_{0})\Subset\Omega_{T}$,
then from \cite[Formulae (4.22) and (4.24)]{AmbPass} one can easily
deduce 
\begin{equation}
\sup_{t\,\in\,(t_{0}-R_{0}^{2},\,t_{0})}\Vert u_{\varepsilon}(\cdot,t)-u(\cdot,t)\Vert_{L^{2}(B_{R_{0}}(x_{0}))}^{2}\rightarrow0\,\,\,\,\,\,\,\,\,\,\mathrm{as}\,\,\varepsilon\rightarrow0^{+},\label{eq:convergenza}
\end{equation}
that is
\[
u_{\varepsilon}\rightarrow u\,\,\,\,\,\,\,\,\mathrm{in}\,\,\,L^{\infty}(t_{0}-R_{0}^{2},t_{0};L^{2}(B_{R_{0}}(x_{0})))\,\,\,\,\,\,\,\,\,\,\mathrm{as}\,\,\varepsilon\rightarrow0^{+}.
\]
\\
Hence, we can conclude that there exists a sequence $\{\varepsilon_{j}\}_{j\in\mathbb{N}}$
such that:
\[
\begin{array}{c}
\bullet\,\,\,0<\varepsilon_{j}\leq1\,\,\,\mathrm{for\,\,\,every\,\,\,}j\in\mathbb{N}\,\,\,\mathrm{and}\,\,\,\varepsilon_{j}\searrow0\,\,\,\mathrm{monotonically\,\,\,as}\,\,\,j\rightarrow+\infty;\\
\\
\bullet\,\,\,u_{\varepsilon_{j}}(x,t)\rightarrow u(x,t)\mathrm{\,\,\,almost\,\,\,everywhere\,\,\,in\,\,\,}Q_{R_{0}}(z_{0})\,\,\,\mathrm{as}\,\,\,j\rightarrow+\infty.\,\,\,\,\,\,\,\,\,\,\,\,\,\,
\end{array}
\]

\vspace{0.3cm}
\section{Physics-informed methodology}\label{sec:method}

\begin{figure}
    \centering
    \includegraphics[width=\textwidth]{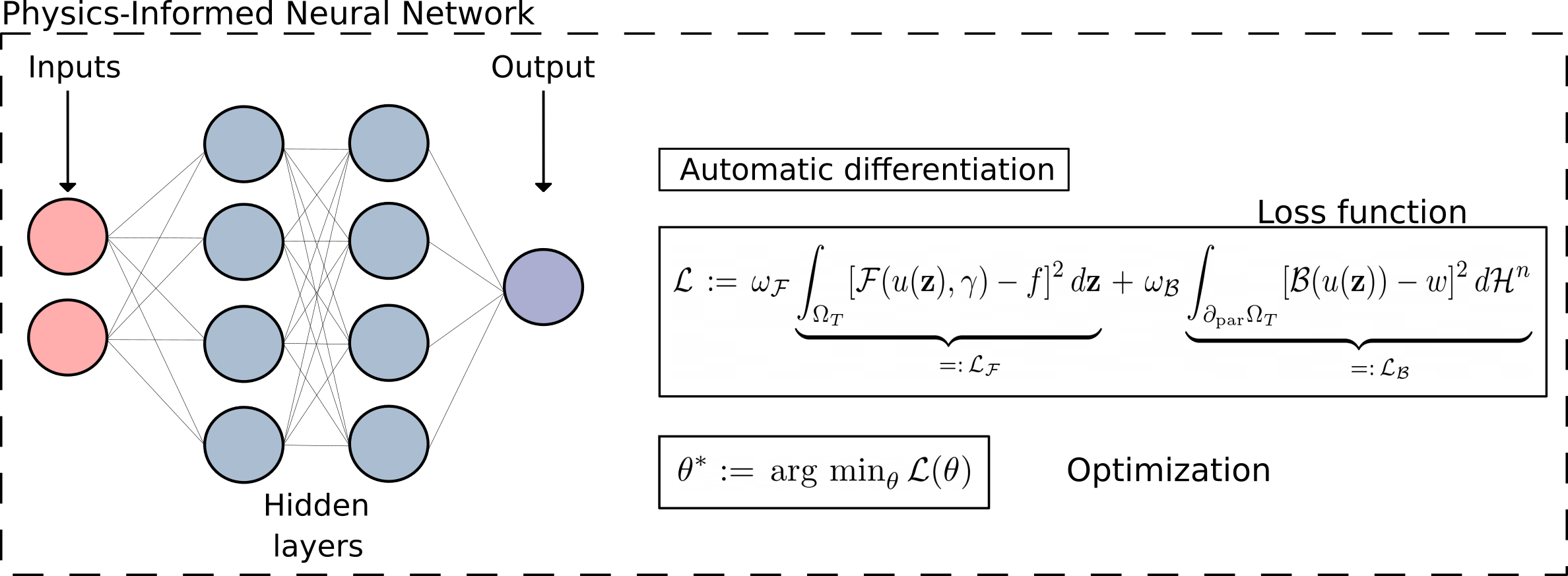}
    \caption{Overall structure of the proposed methodology. An FF-DNN serves as the neural network's architecture. Automatic differentiation is employed to calculate the loss terms associated with the model's dynamics. The loss function comprises two components: the physics loss, represented by $\mathcal{L}_{\mathcal{F}}$, and the boundary loss denoted by $\mathcal{L}_{\mathcal{B}}$. During the optimization phase, the objective is to minimize the loss function with respect to the set of hyperparameters $\theta$.}
    \label{pipeline}
\end{figure}

\noindent $\hspace*{1em}$PINNs
are a type of SciML approach used in neural
networks to solve PDEs. Unlike traditional
neural networks, PINNs incorporate physics constraints into the model,
resulting in predicted outcomes that adhere more closely to the natural
laws governing the specific problem being addressed. The general form
of the problem involves a PDE along with initial and/or boundary conditions.\\
$\hspace*{1em}$In particular, we consider a (well-posed) problem
of the type
\begin{equation}
\begin{cases}
\mathcal{F}(u(\mathbf{x},t),\gamma)=f\hspace{0.5cm} & \text{ if }(\mathbf{x},t)\in\Omega_{T}:=\Omega\times(0,T)\\
\mathcal{B}(u(\mathbf{x},t))=w & \text{ if }(\mathbf{x},t)\in\partial_{\mathrm{par}}\Omega_{T},
\end{cases}\label{differential problem}
\end{equation}
where $\Omega$ is a bounded domain in $\mathbb{R}^{n}$, $\mathcal{F}$
denotes a nonlinear differential operator, $\gamma$ is a parameter
associated with the physics of the problem, $\mathcal{B}$ is an operator
defining arbitrary initial-boundary conditions, the functions $f$
and $w$ represent the problem data, while $u(\mathbf{x},t)$ denotes
the unknown solution.\\
$\hspace*{1em}$The objective of PINNs
is to predict the solution to (\ref{differential problem}) by training
the neural network to minimize a cost function. The neural network's
architecture used for PINNs is typically a FeedForward fully-connected Neural Network
(FF-DNN), also known as Multi-Layer Perceptron (MLP). In an FF-DNN,
information flows only forward direction, in the sense that the neural network does not form a loop.
Furthermore, all neurons are interconnected. Once the number $N$ of
hidden layers has been chosen, for any $i\in\{1,\ldots,N\}$ and set
$\mathbf{z}=(\mathbf{x},t)$ we define 
\[
\Gamma_{i}(\mathbf{z}_{i-1}):=\,W_{i}\,\mathbf{z}_{i-1}+\mathbf{b}_{i}\,,
\]
where $W_{i}$ is the weights matrix of the links between the layers
$i-1$ and $i$, while $\mathbf{b}_{i}$ corresponds to the biases
vector. Then, a generic layer of the neural network is defined by
\[
h_{i}(\mathbf{z}_{i-1};W_{i},\mathbf{b}_{i}):=\,\varphi_{i}(\Gamma_{i}(\mathbf{z}_{i-1})),\,\,\,\,\,\,\,\,\,\,\,i\in\{1,\ldots,N\},
\]
for some nonlinear activation function $\varphi_{i}$. The output
of the FF-DNN, denoted by $\hat{u}_{\theta}(\mathbf{z})$, can be
expressed as a composition of these layers by 
\begin{equation}
\hat{u}_{\theta}(\mathbf{\mathbf{z}}):=\,(\Gamma_{N}\circ\varphi\circ\Gamma_{N-1}\circ\cdots\circ\varphi\circ\Gamma_{1})(\mathbf{z})\,,
\end{equation}
where $\theta$ represents the set of hyperparameters of the neural network and the activation function $\varphi$ is assumed to be the same
for all layers. To solve the differential problem (\ref{differential problem})
using PINNs, the PDE is approximated by finding an optimal set $\theta^{*}$
of neural network hyperparameters that minimizes a loss function $\mathcal{L}$.
This function consists of two components: the former, denoted by $\mathcal{L}_{\mathcal{F}}$,
is related to the differential equation, while the latter, here denoted
by $\mathcal{L}_{\mathcal{B}}$, is connected to the initial-boundary
conditions (see Fig. \ref{pipeline}). In particular, the loss function can be defined as follows
\begin{equation}
\mathcal{L}\,:=\,\omega_{\mathcal{F}}\underbrace{\int_{\Omega_{T}}[\mathcal{F}(u(\mathbf{z}),\gamma)-f]^{2}\,d\mathbf{z}}_{\text{\,\,\,\,\,\,\,\,\,\,\,=:\,\ensuremath{\mathcal{L}_{\mathcal{F}}}}}\,+\,\,\omega_{\mathcal{B}}\underbrace{\int_{\partial_{\mathrm{par}}\Omega_{T}}[\mathcal{B}(u(\mathbf{z}))-w]^{2}\,d\mathcal{H}^{n}}_{\text{\,\,\,\,\,\,\,\,\,\,\,=:\,\ensuremath{\mathcal{L}_{\mathcal{B}}}}},\label{loss function}
\end{equation}
where $\omega_{\mathcal{F}}$ and $\omega_{\mathcal{B}}$ represent the weights that are usually
applied to balance the importance of each component. Hence, we can
write 
\begin{equation}
\theta^{*}:=\,\argmin_{\theta} \,\mathcal{L}(\theta).
\end{equation}
The aim of this approach is to approximate the solution
of the PDE satisfying the initial-boundary conditions. This is known in the literature as the \textit{direct
problem}, which is the only one we will address here.\vspace{0.6cm}

\section{Numerical results}\label{sec:results}

\noindent $\hspace*{1em}$In this section we evaluate the accuracy
and effectiveness of our predictive method, by testing it with five
problems of the type (\ref{eq:mainprob}) whose exact solutions are
known. For each problem, we will denote the \textit{exact solution}
by $u$, and the \textit{predicted} (or \textit{approximate}) \textit{solution}
by $\hat{u}$. Sometimes, by abuse of language, for a given time $t\geq0$
we will refer to the partial maps $u(\cdot,t)$ and $\hat{u}(\cdot,t)$
as the exact solution and the predicted (or approximate) solution
respectively. The meaning will be clear from the context every time.
We will deal with each test problem separately, so that no confusion
can arise. In the first three problems, $\Omega$ will be a bounded
domain of $\mathbb{R}^{2}$, while, in the last two problems, $\Omega$
will denote the open unit sphere of $\mathbb{R}^{3}$ centered at
the origin.\\
$\hspace*{1em}$ In addition, for each of the test problems, we employed
the same neural network architecture. This consists of four layers,
each with $20$ neurons. We utilized the hyperbolic tangent function
as the activation function for both the input layer and the hidden
layers, while a linear function served as the activation function
for the output layer. Lastly, to train the neural network, we conducted
$80000$ epochs with a learning rate (lr) of $3\times10^{-3}$ and employed
the Adaptive Moment Estimation (ADAM) optimizer. The decision to set the lr to the constant value $3\times10^{-3}$ was based on the observation that this specific hyperparameter led to the optimal convergence of our method. Experimentation with lr set to $1\times10^{-1}$ highlighted the network's inability to achieve convergence, while using an lr of $1\times10^{-5}$ allowed the method to converge, albeit requiring a significantly higher number of epochs. The latter scenario, while ensuring convergence, proved to be less computationally efficient. The experiments were performed on a NVIDIA GeForce
RTX $3080$ GPU with AMD Ryzen $9$ $5950$X $16$-Core Processor
and $128$ GB of RAM.
\subsection{First test problem}
\begin{figure}[H]
\centering{}\includegraphics[scale=0.85]{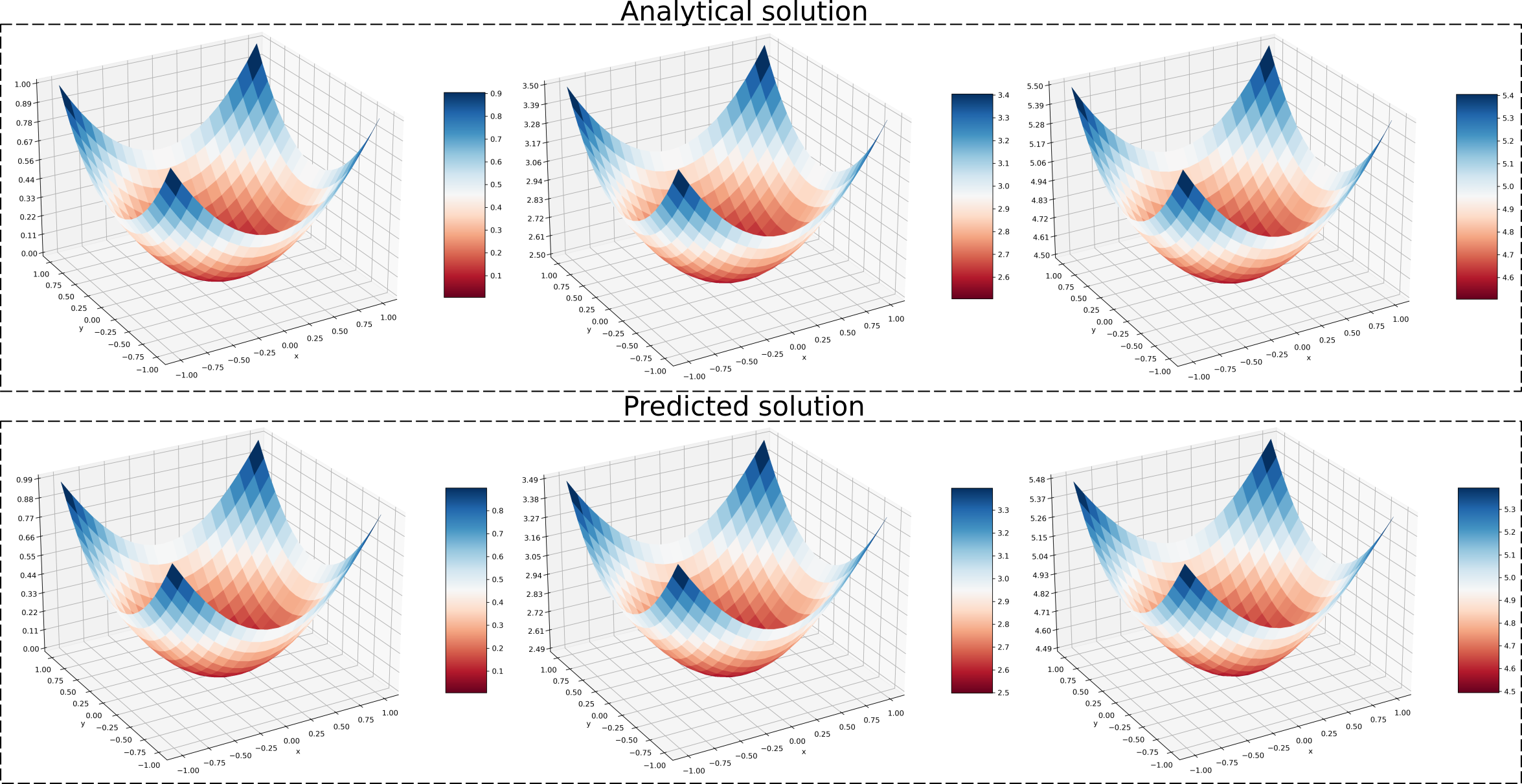}\caption{Plots of the exact solution to problem (\ref{eq:P1}) (above) and
the predicted solution $\hat{u}(\cdot,t)$ (below) for $t=0$ (left),
$t=2.5$ (center) and $t=4.5$ (right).}
\label{imm1}
\end{figure}

\begin{figure}[H]
\centering{}\includegraphics[scale=0.47]{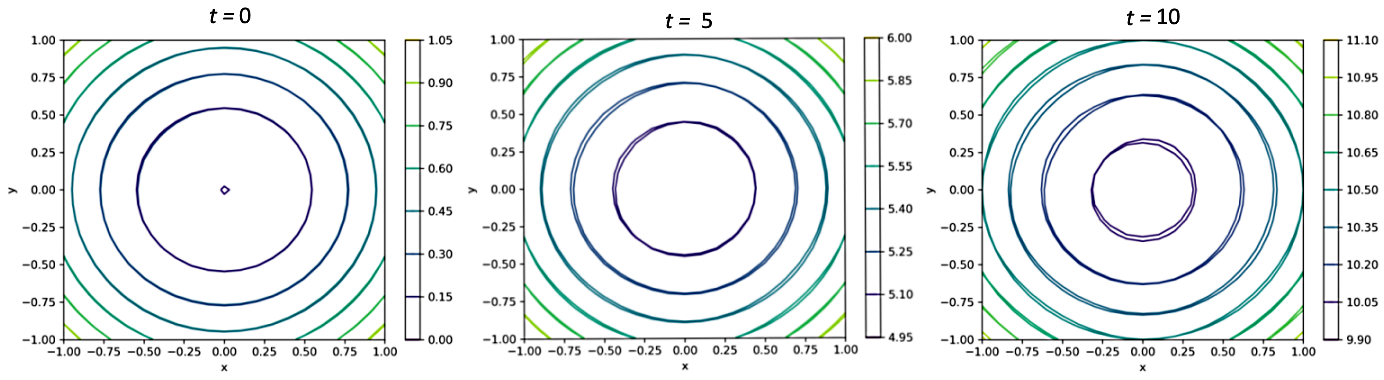}\caption{Superposition of the level curves of the exact solution $u(\cdot,t)$
and the predicted solution $\hat{u}(\cdot,t)$ for $t=0$ (left),
$t=5$ (center) and $t=10$ (right). The contour lines corresponding
to the same level are almost indistinguishable for any fixed $t\in[0,10]$.}
\label{imm2}
\end{figure}

\noindent $\hspace*{1em}$The first test problem that we consider is\\
\begin{equation}
\begin{cases}
\begin{array}{cc}
\partial_{t}v-\mathrm{div}\left((\vert\nabla v\vert-1)_{+}\,\frac{\nabla v}{\vert\nabla v\vert}\right)=1 & \mathrm{in\,\,}\Omega_{T},\,\,\,\,\,\,\,\,\,\,\,\,\,\,\,\,\,\,\,\,\,\,\,\,\,\,\,\,\,\,\,\,\,\,\,\,\,\,\,\,\,\,\,\,\,\,\,\,\,\\
v(x,y,0)=\frac{1}{2}(x^{2}+y^{2}) & \mathrm{if}\,\,(x,y)\in\overline{\Omega},\,\,\,\,\,\,\,\,\,\,\,\,\,\,\,\,\,\,\,\,\,\,\,\,\,\,\,\,\,\,\,\,\\
v(x,y,t)=\frac{1}{2}+t\,\,\,\,\,\,\,\,\,\,\,\,\, & \,\,\,\mathrm{if}\,\,(x,y)\in\partial\Omega\,\,\mathrm{\land}\,\,t\in(0,T),
\end{array} & \tag{P1}\end{cases}\label{eq:P1}
\end{equation}
\\
where $\Omega=\{(x,y)\in\mathbb{R}^{2}:x^{2}+y^{2}<1\}$. The exact
solution of this problem is given by 
\[
u(x,y,t)=\frac{1}{2}(x^{2}+y^{2})+t.
\]
Therefore, for any fixed time $t\geq0$ the graph of the function
$u(\cdot,t)$ is an elliptic paraboloid. As time goes on, this paraboloid
slides along an oriented vertical axis at a constant velocity, without
deformation, since $\partial_{t}u\equiv1$ over $\Omega_{T}$ (see
Fig. \ref{imm1}, above).\\
$\hspace*{1em}$To train the neural network,
in each experiment we have initially used $441$ points to suitably
discretize the domain $\Omega$ and its boundary $\partial\Omega$,
and $21$ equispaced points in the time interval $[0,T]$. Once the
network has been trained, we have made a prediction of the solution
to problem (\ref{eq:P1}) at different times $t$ (Fig. \ref{imm1},
below).
\noindent $\hspace*{1em}$What has been verified is that the plot
of the predicted solution $\hat{u}(\cdot,t$) has precisely the same
shape and geometric properties as the graph of the exact solution $u(\cdot,t)$, for both short and long
times $t$. Moreover, the time evolution of the approximate solution
$\hat{u}$ exactly mirrors the behavior described for the known solution
$u$. A further interesting aspect that can be noticed is that the
level curves of the approximate solution $\hat{u}(\cdot,t$) overlap
almost perfectly those of $u(\cdot,t$), provided that $t$ is not
very large (see Fig. \ref{imm2}).
\vspace{5mm}
\newpage
\noindent $\hspace*{1em}$We have also noted that, at time $t=0$,
the approximate solution is basically equal to zero in a very tiny
region around the origin $(0,0)$ of the $xy$-plane. This means that
the said region is composed of ``numerical zeros'' of the solution
predicted at time $t=0$, while we know that $u(x,y,0)=0$ if and
only if $(x,y)=(0,0)$. However, this discrepancy is actually negligible,
since the order of magnitude of $u(x,y,0)$ is not greater than $10^{-6}$
within the above region.\\
$\hspace*{1em}$To assess the accuracy of our predictive method and
the numerical convergence of the solution $\hat{u}$ toward $u$ in
a more quantitative way, we now look at the time behavior of the $L^{2}$-error
$\|\hat{u}(\cdot,t)-u(\cdot,t)\|_{L^{2}(\Omega)}$ by considering
the natural quantities\\
\\

\noindent 
\begin{equation}
\mathcal{E}(T):=\underset{t\,\in\,(0,T)}{\sup}\,\|\hat{u}(\cdot,t)-u(\cdot,t)\|_{L^{2}(\Omega)}^{2}\label{eq:L2err}
\end{equation}
and\\
\begin{equation}
\mathcal{E}_{rel}(T):=\,\frac{\mathcal{E}(T)}{\underset{t\,\in\,(0,T)}{\sup}\,\|u(\cdot,t)\|_{L^{2}(\Omega)}^{2}}\,.\label{eq:err_rel}
\end{equation}
\\
Passing from Cartesian to polar coordinates,
one can easily find that\\
\\
\[
\|u(\cdot,t)\|_{L^{2}(\Omega)}^{2}\,=\underset{\Omega\,\,\,\,}{\iint}\left[\frac{1}{2}(x^{2}+y^{2})+t\right]^{2}dx\,dy\,=\,\pi\,\left(t^{2}+\frac{t}{2}+\frac{1}{12}\right),
\]
\\
and therefore
\[
\mathcal{E}_{rel}(T)=\,\frac{12\cdot\mathcal{E}(T)}{\pi\,(12\,T^{2}+6\,T+1)}\,.
\]
\\
\\
During our numerical experiments, we have estimated both $\mathcal{E}(T)$
and $\mathcal{E}_{rel}(T)$ for 
\[
T\in\{1,10,100,200,300,400,500\}.
\]
The results that we have obtained are shown in Table \ref{tab1}.
The estimates of $\mathcal{E}(1)$ and $\mathcal{E}(10)$ are equal
to $2.24\times10^{-5}$ and $4.30\times10^{-4}$ respectively, which
is very satisfactory, especially considering that the order of magnitude
of $\mathcal{E}_{rel}(1)$ and $\mathcal{E}_{rel}(10)$ is equal to
$10^{-6}$.\\
In order to get more accurate estimates for larger values
of $T$, for every fixed $T\geq100$ we have used
$2.1\times T$ equispaced points (instead of the initial $21$) to
discretize the time interval $[0,T]$. By doing so, we have observed
that the variation of the estimate of $\mathcal{E}(T)$ displays a
monotonically increasing behavior, in accordance with the definition
(\ref{eq:L2err}). However, even for $100\leq T\leq500$,
the order of magnitude of $\mathcal{E}_{rel}(T)$ remains not greater
than $10^{-6}$. Therefore, for this first test problem, we can conclude
that our predictive method is indeed very accurate and efficient,
on both a short and long time scale. 
\newpage
\begin{table}[H]
\noindent \centering{}%
\begin{tabular}{|c|c|c|}
\hline 
\textbf{Final time $\boldsymbol{T}$} & \textbf{Estimate of} $\boldsymbol{\mathcal{E}(T)}$ & \textbf{Estimate of} $\boldsymbol{\mathcal{E}_{rel}(T)}$\tabularnewline
\hline 
$1$ & $2.24\times10^{-5}$ & $4.50\times10^{-6}$\tabularnewline
\hline 
$10$ & $4.30\times10^{-4}$ & $1.30\times10^{-6}$\tabularnewline
\hline 
$100$ & $1.87\times10^{-1}$ & $5.92\times10^{-6}$\tabularnewline
\hline 
$200$ & $2.05\times10^{-1}$ & $1.63\times10^{-6}$\tabularnewline
\hline 
$300$ & $2.45\times10^{-1}$ & $8.65\times10^{-7}$\tabularnewline
\hline 
$400$ & $2.47\times10^{-1}$ & $4.91\times10^{-7}$\tabularnewline
\hline 
$500$ & $2.80\times10^{-1}$ & $3.56\times10^{-7}$\tabularnewline
\hline 
\end{tabular}\caption{Estimates of $\mathcal{E}(T)$ and $\mathcal{E}_{rel}(T)$
for $T\in\{1,10,100,200,300,400,500\}$.}
\label{tab1}
\end{table}

\begin{figure}[H]
\centering{}\includegraphics[scale=0.82]{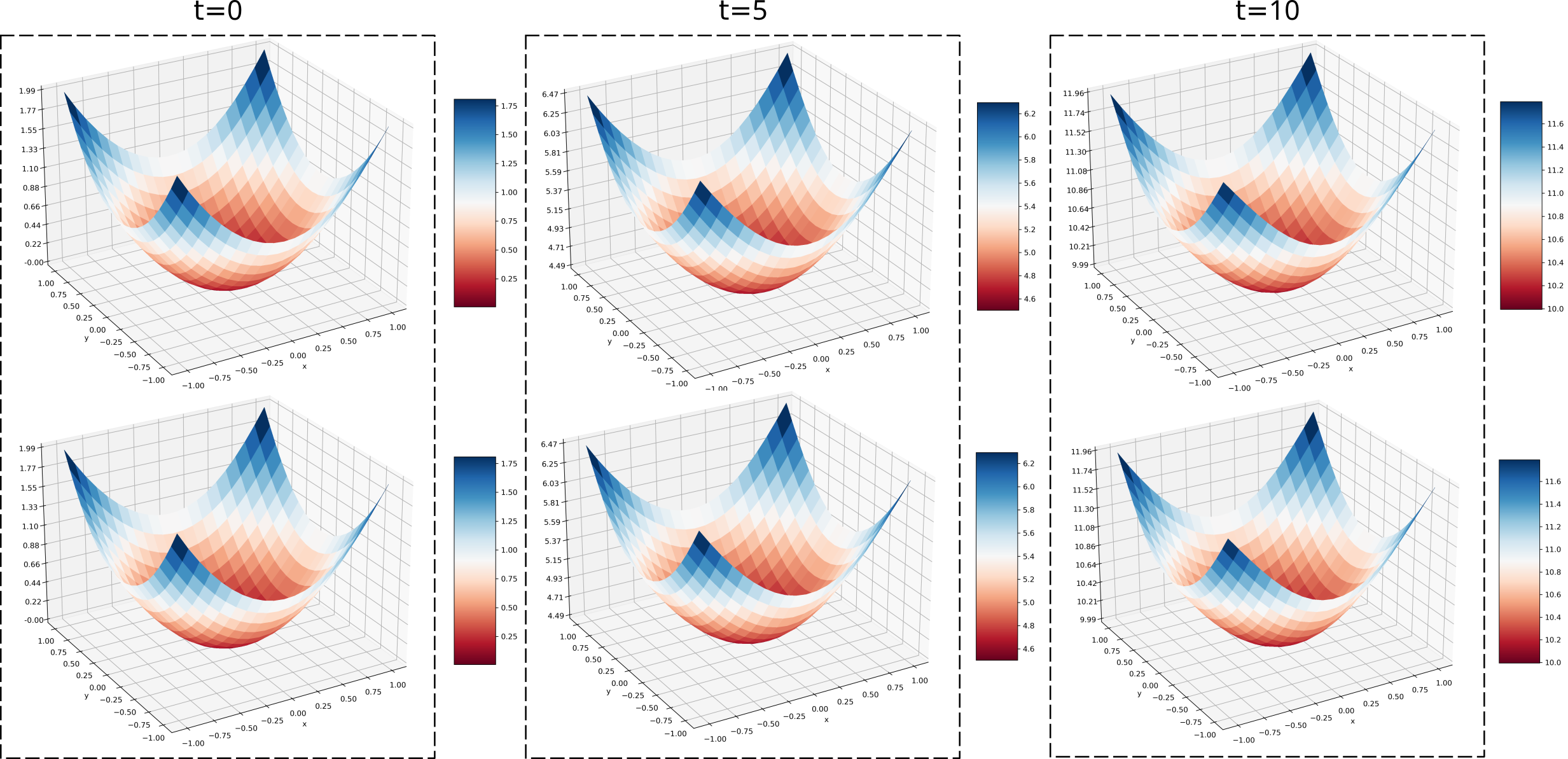}\caption{Plots of the approximate solution $\hat{u}_{\varepsilon}(\cdot,t)$
for $\varepsilon=10^{-3}$ (above) and $\varepsilon=10^{-9}$
(below), at times $t=0$ (left), $t=5$ (center) and $t=10$ (right).
Here $\Omega'=\Omega$ and $[t_{1},t_{2}]=[0,10]$.}
 \label{app2}
\end{figure}

\begin{figure}[H]
\centering{}\includegraphics[scale=0.93]{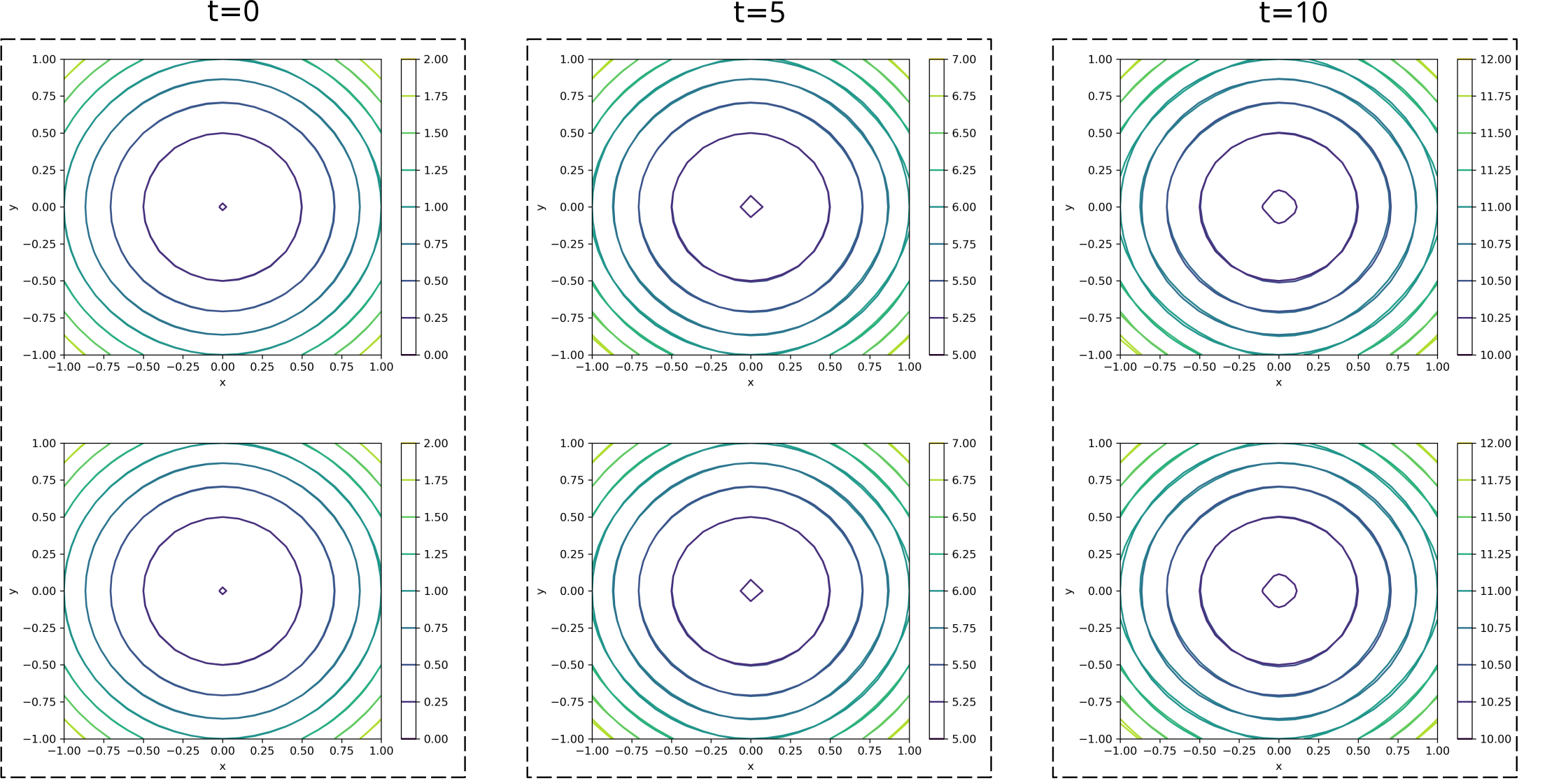}\caption{Superposition of the level curves of the exact solution $u(\cdot,t)$
and the predicted solution $\hat{u}_{\varepsilon}(\cdot,t)$ for $\varepsilon=10^{-3}$
(above) and $\varepsilon=10^{-9}$ (below), at times $t=0$ (left),
$t=5$ (center) and $t=10$ (right). Here $\Omega'=\Omega$ and $[t_{1},t_{2}]=[0,10]$.
For every $\varepsilon\in[10^{-9},10^{-3}]$, the contour lines of
the approximate solution $\hat{u}_{\varepsilon}(\cdot,t)$ perfectly
overlap those of $u(\cdot,t)$ for any fixed $t\in[0,10]$.}
\label{app1}
\end{figure}

\noindent $\hspace*{1em}$Now, for $0<\varepsilon\leq1$ we consider
the problem
\\
\begin{equation}
\begin{cases}
\begin{array}{cc}
\partial_{t}v_{\varepsilon}-\mathrm{div}\left(H(\nabla v_{\varepsilon})+\varepsilon\,\nabla v_{\varepsilon}\right)=1 & \mathrm{in\,\,}\Omega'_{1,2}:=\Omega'\times(t_{1},t_{2}),\,\,\,\,\,\,\,\,\,\\
v_{\varepsilon}(x,y,t_{1})=\frac{1}{2}(x^{2}+y^{2})+t_{1} & \mathrm{if}\,\,(x,y)\in\overline{\Omega'},\,\,\,\,\,\,\,\,\,\,\,\,\,\,\,\,\,\,\,\,\,\,\,\,\,\,\,\,\,\,\,\,\\
v_{\varepsilon}(x,y,t)=\frac{1}{2}(x^{2}+y^{2})+t & \,\,\,\,\,\mathrm{if}\,\,(x,y)\in\partial\Omega'\,\,\mathrm{\land}\,\,t\in(t_{1},t_{2}),
\end{array}\end{cases}\label{eq:appro1}
\end{equation}
\\
which is nothing but the approximating problem
(\ref{eq:CAUCHYDIR}) associated with (\ref{eq:P1}). In what follows,
we will denote the exact solution of (\ref{eq:appro1}) by $u_{\varepsilon}$,
while the predicted solution will be denoted by $\hat{u}_{\varepsilon}$.\\
Throughout our tests, for $10^{-9}\leq\varepsilon\leq10^{-3},$
for $\Omega'=\Omega$ and for $(t_{1},t_{2})=(0,T)$, we have observed
that the plots of the predicted solution $\hat{u}_{\varepsilon}(\cdot,t$) and the exact
solution $u(\cdot,t)$ share the same configurations and geometric peculiarities, on both a short
and long time scale (see, e.g. Figure \ref{app2}).
Furthermore, we have seen that the evolution over time of $\hat{u}_{\varepsilon}$ reflects the behavior depicted
for the solution $u$ quite faithfully. In addition, the contour lines
of $\hat{u}_{\varepsilon}(\cdot,t)$ perfectly overlap those of $u(\cdot,t)$,
at least for not very long times $t$ (see Fig. \ref{app1}).

\noindent $\hspace*{1em}$Let us now assume that
\[
\Omega'_{1,2}=Q_{\frac{1}{2}}(z_{0}):=B_{\frac{1}{2}}(\mathbf{0})\times\left(\frac{7}{4},2\right),
\]
where $z_{0}=(\mathbf{0},2)=(0,0,2)$. Then, the limit in (\ref{eq:convergenza})
suggests that $\hat{u}_{\varepsilon}$ should numerically converge
to $u$ as $\varepsilon\searrow0$. To obtain a numerical evidence
of such convergence, we have chosen $\Omega'=B_{1/2}(\mathbf{0})$
and $(t_{1},t_{2})=(\frac{7}{4},2)$ into (\ref{eq:appro1}) and examined the time behavior of the $L^{2}$-error $\|\hat{u}_{\varepsilon}(\cdot,t)-u(\cdot,t)\|_{L^{2}(\Omega')}$,
by evaluating the quantities\\ 
\[
E_{\varepsilon}:=\underset{t\,\in\,(\frac{7}{4},\,2)}{\sup}\,\|\hat{u}_{\varepsilon}(\cdot,t)-u(\cdot,t)\|_{L^{2}(B_{1/2}(\mathbf{0}))}^{2}
\]
and
\[
E_{rel}\,\equiv\,E_{rel}(\varepsilon)\,:=\,\frac{E_{\varepsilon}}{\underset{t\,\in\,(\frac{7}{4},\,2)}{\sup}\,\|u(\cdot,t)\|_{L^{2}(B_{1/2}(\mathbf{0}))}^{2}}\,.
\]\\
Switching from Cartesian to polar coordinates, one can easily compute\\
\[
\|u(\cdot,t)\|_{L^{2}(B_{1/2}(\mathbf{0}))}^{2}\,=\underset{B_{\frac{1}{2}}(\mathbf{0})\,\,\,\,}{\iint}\left[\frac{1}{2}(x^{2}+y^{2})+t\right]^{2}dx\,dy\,=\,\frac{\pi}{4}\left(t^{2}+\frac{t}{8}+\frac{1}{192}\right),
\]
from which it immediately follows that 
\[
E_{rel}\,=\,\frac{768\,E_{\varepsilon}}{817\,\pi}\,.
\]
In the testing phase, we have estimated both $E_{\varepsilon}$
and $E_{rel}$ for $\varepsilon\in\{10^{-15},10^{-14},\ldots,10^{-1},1\}$.
Table \ref{tab2} shows the results obtained and reveals that the predicted solution $\hat{u}_{\varepsilon}$ converges to $u$ as $\varepsilon$ tends to zero, although not very quickly. In fact, the estimates of  $E_{\varepsilon}$
and $E_{rel}(\varepsilon)$ approach zero with a convergence rate much lower than that of $\varepsilon$. Furthermore, they seem to start decreasing monotonically, i.e. without oscillations, for $\varepsilon \leq 10^{-10}$.
\vspace{3mm}
\begin{table}[H]
\noindent \centering{}%
\begin{tabular}{|c|c|c|}
\hline 
\textbf{$\boldsymbol{\varepsilon}$} & \textbf{Estimate of} $\boldsymbol{E_{\varepsilon}}$ & \textbf{Estimate of} $\boldsymbol{E}_{\boldsymbol{rel}}$
\tabularnewline
\hline 
$1$ & $5.816\times10^{-4}$ & $1.74\times10^{-4}$
\tabularnewline
\hline 
$10^{-1}$ & $1.658\times10^{-5}$ & $4.96\times10^{-6}$
\tabularnewline
\hline 
$10^{-2}$ & $2.485\times10^{-5}$ & $7.44\times10^{-6}$
\tabularnewline
\hline 
$10^{-3}$ & $9.280\times10^{-6}$ & $2.78\times10^{-6}$\tabularnewline
\hline 
$10^{-4}$ & $5.331\times10^{-6}$ & $1.60\times10^{-6}$\tabularnewline
\hline 
$10^{-5}$ & $9.681\times10^{-6}$ & $2.90\times10^{-6}$\tabularnewline
\hline 
$10^{-6}$ & $7.175\times10^{-6}$ & $2.15\times10^{-6}$\tabularnewline
\hline 
$10^{-7}$ & $2.027\times10^{-6}$ & $6.06\times10^{-7}$\tabularnewline
\hline 
$10^{-8}$ & $8.205\times10^{-6}$ & $2.45\times10^{-6}$\tabularnewline
\hline 
$10^{-9}$ & $3.749\times10^{-6}$ & $1.12\times10^{-6}$
\tabularnewline
\hline 
$10^{-10}$ & $4.254\times10^{-6}$ & $1.27\times10^{-6}$
\tabularnewline
\hline 
$10^{-11}$ & $3.611\times10^{-6}$ & $1.08\times10^{-6}$
\tabularnewline
\hline 
$10^{-12}$ & $2.769\times10^{-6}$ & $8.28\times10^{-7}$
\tabularnewline
\hline 
$10^{-13}$ & $2.399\times10^{-6}$ & $7.18\times10^{-7}$
\tabularnewline
\hline 
$10^{-14}$ & $2.017\times10^{-6}$ & $6.03\times10^{-7}$
\tabularnewline
\hline 
$10^{-15}$ & $1.649\times10^{-6}$ & $4.94\times10^{-7}$
\tabularnewline
\hline 
\end{tabular}\caption{Estimates of $E_{\varepsilon}$ and $E_{rel}(\varepsilon$) for $\varepsilon\in\{10^{-15},10^{-14},\ldots,10^{-1},1\}$.}
\label{tab2}
\end{table}

\subsection{Second test problem}
\begin{figure}[H]
\centering{}\includegraphics[scale=0.85]{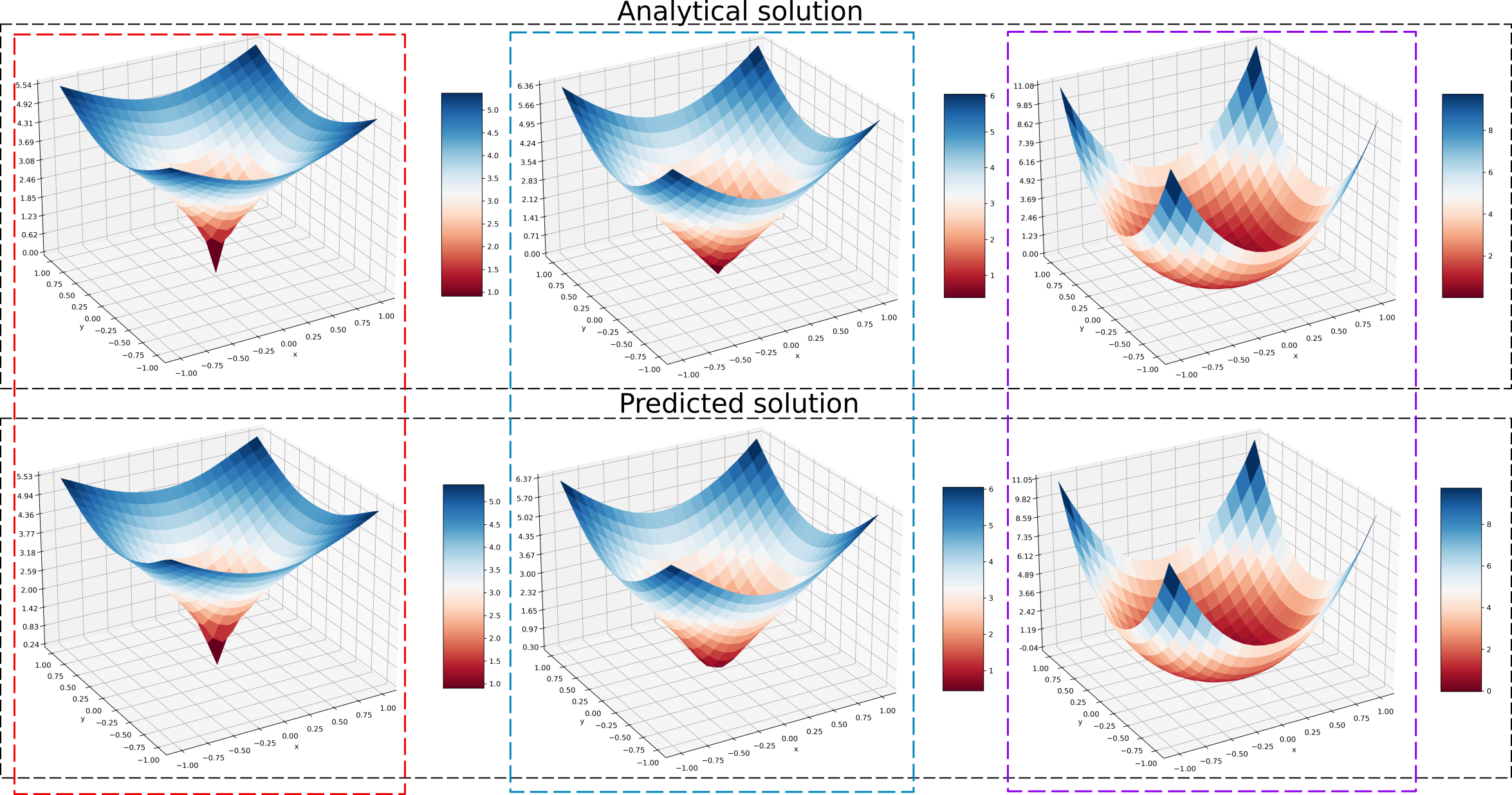}\caption{Plots of the exact solution (above) and the predicted solution (below)
to problem (\ref{eq:P2}) at time $t=4.5$ for $\textcolor{red}{\alpha=0.3}$
(left), $\textcolor[RGB]{0,134,185}{\alpha=0.5}$ (center) and $\textcolor[RGB]{141,0,237}{\alpha=1.3}$
(right).}
\label{imm2_1}
\end{figure}

\noindent $\hspace*{1em}$Let $\alpha>0$. As a second test problem
we consider\\
\\
\begin{equation}
\begin{cases}
\begin{array}{cc}
\partial_{t}v-\mathrm{div}\left((\vert\nabla v\vert-1)_{+}\,\frac{\nabla v}{\vert\nabla v\vert}\right)=f & \mathrm{in\,\,}\Omega_{T},\,\,\,\,\,\,\,\,\,\,\,\,\,\,\,\,\,\,\,\,\,\,\,\,\,\,\,\,\,\,\,\,\,\,\,\,\,\,\,\,\,\,\,\,\,\,\,\,\,\\
v(x,y,0)=0 & \mathrm{if}\,\,(x,y)\in\overline{\Omega},\,\,\,\,\,\,\,\,\,\,\,\,\,\,\,\,\,\,\,\,\,\,\,\,\,\,\,\,\,\,\,\,\\
v(x,y,t)=t & \,\,\,\,\mathrm{if}\,\,(x,y)\in\partial\Omega\,\,\mathrm{\land}\,\,t\in(0,T),
\end{array} & \tag{P2}\end{cases}\label{eq:P2}
\end{equation}
\\
\\
where $\Omega=\{(x,y)\in\mathbb{R}^{2}:x^{2}+y^{2}<1\}$
again and
\vspace{5mm}
\[
{\normalcolor {\color{blue}{\normalcolor f(x,y,t):=\begin{cases}
\begin{array}{cc}
(x^{2}+y^{2})^{\alpha} & \,\,\mathrm{if}\,\,2\alpha\,t\,(x^{2}+y^{2})^{\alpha-\frac{1}{2}}\leq1,\\
(x^{2}+y^{2})^{\alpha}-4\alpha^{2}t\,(x^{2}+y^{2})^{\alpha-1}+\frac{1}{\sqrt{x^{2}+y^{2}}} & \,\,\mathrm{if}\,\,2\alpha\,t\,(x^{2}+y^{2})^{\alpha-\frac{1}{2}}>1.
\end{array} & {\normalcolor }\end{cases}}}}
\]
\\
\\
The exact solution of problem (\ref{eq:P2})
is given by
\[
u(x,y,t)\equiv u_{\alpha}(x,y,t):=t\,(x^{2}+y^{2})^{\alpha}.
\]
At any fixed time $t>0$, the shape and geometric properties of the
graph of $u(\cdot,t)$ strongly depend on the value of the parameter
$\alpha$.\\
$\hspace*{1em}$If $\alpha=\frac{1}{2}$,
then the graph of $u(\cdot,t)$ is a cone whose vertex coincides with
the origin $(0,0,0)$ at any given positive time $t$. As time goes
on, the cone in question gets narrower and narrower around the vertical
axis. In this case, the plot of the approximate solution $\hat{u}(\cdot,t)$
has the same form as the graph of the exact solution $u(\cdot,t)$
for both short and long times $t>0$, except
near the origin, where the tip of the cone appears to have been smoothed
out (see Figure \ref{imm2_1}, center). However, this is not a surprise
at all, since we already know that for $t>0$ the function 
\[
(x,y)\in\Omega\mapsto t\,(x^{2}+y^{2})^{\frac{1}{2}}
\]
is not differentiable at the center $(0,0)$ of $\Omega$.

\noindent 

\begin{figure}[H]
\centering{}\includegraphics[scale=0.83]{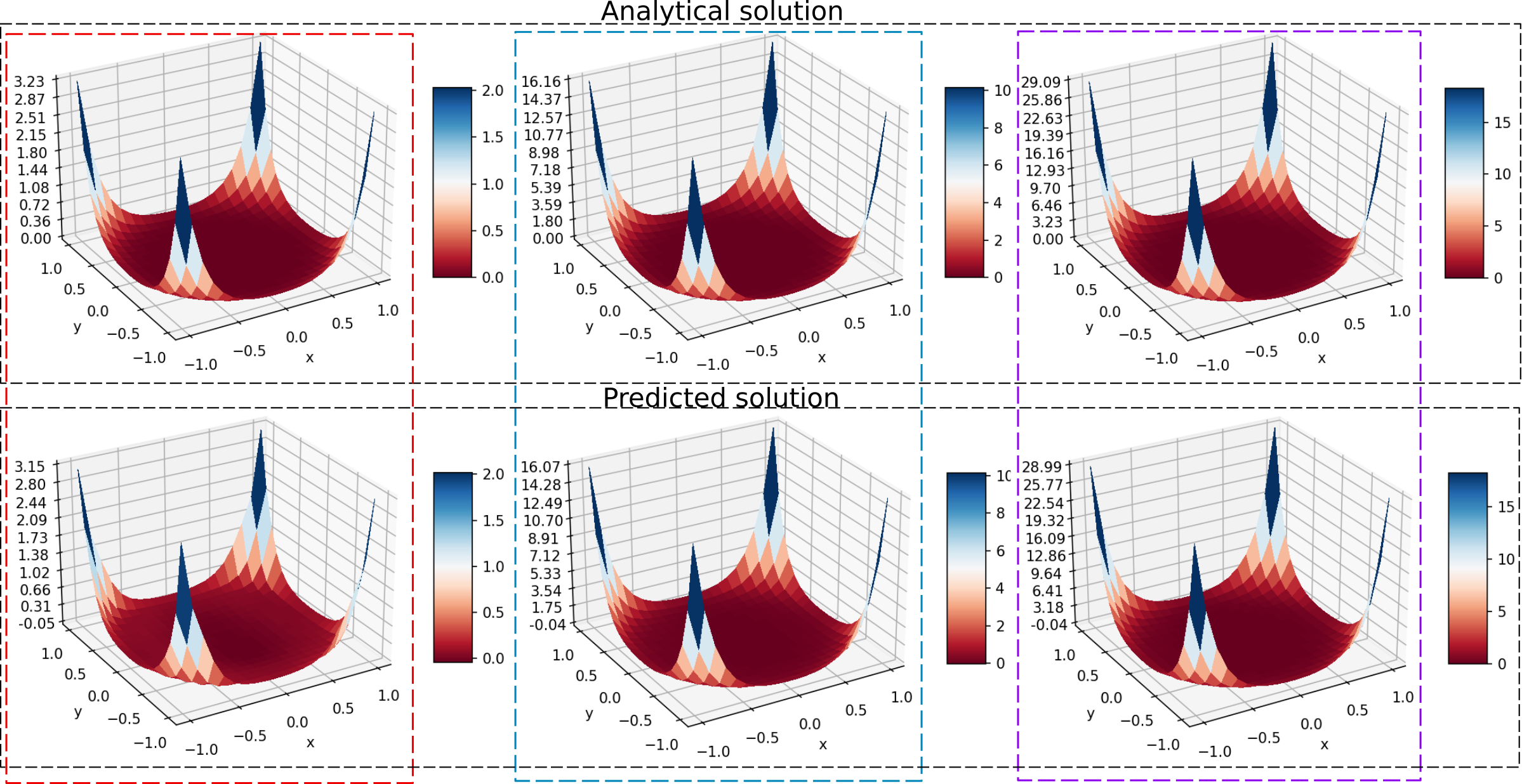}\caption{Plots of the exact solution to problem (\ref{eq:P2}) (above) and
the predicted solution $\hat{u}(\cdot,t)$ (below) for $\alpha=5$,
at times $\textcolor{red}{t=0.5}$ (left), $\textcolor[RGB]{0,134,185}{t=5}$
(center) and $\textcolor[RGB]{141,0,237}{t=10}$ (right).}
\label{imm2_2}
\end{figure}

\noindent $\hspace*{1em}$When $0<\alpha<\frac{1}{2}$, the graph
of $u(\cdot,t)$ is cusp-shaped for any fixed time $t>0$, the origin
now being a cusp for all positive times. In this case, a loss on convexity
occurs, which is also observed in the plot of the predicted solution $\hat{u}(\cdot,t$) for all times $t>0$
(see, e.g. Figure \ref{imm2_1}, left).\\
$\hspace*{1em}$Lastly, when $\alpha>\frac{1}{2}$ the graph of $u(\cdot,t)$
is no longer cusp-shaped and becomes increasingly narrow around the
vertical axis as $t$ increases. Furthermore, for any fixed $t>0$
the exact solution $u(\cdot,t)$ is convex again and its graph gets
flatter and flatter near the origin when $\alpha>>1$ (see Figure
\ref{imm2_2}).\\
$\hspace*{1em}$In all three of the above cases, we have noticed that the plot of $\hat{u}(\cdot,t$) is basically identical in its shape and geometry to the graph of the exact solution $u(\cdot,t)$,
for both short and long periods $t$.\\
Moreover, also for problem (\ref{eq:P2}) we have verified that the
time evolution of the predicted solution
faithfully reflects the trend described for the exact solution
in all three previous cases. Therefore, we may conclude that $\alpha=\frac{1}{2}$
represents a critical value for the global behavior of both the exact
and the predicted solution.

\begin{figure}[H]
\centering{}\includegraphics[scale=1.04]{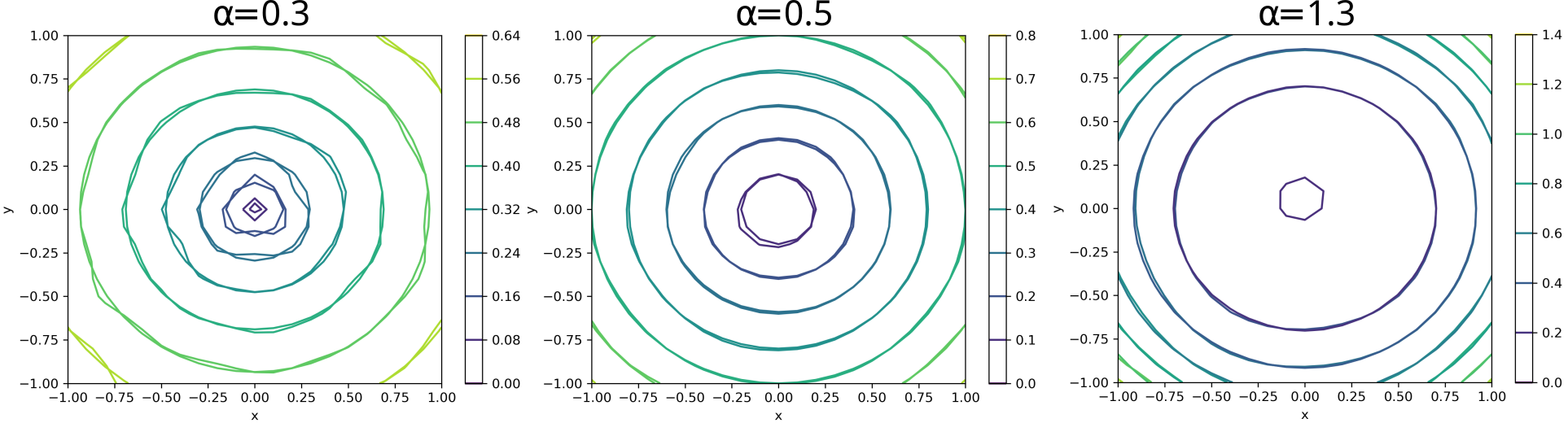}\caption{Superposition of the level curves of the exact solution $u_{\alpha}(\cdot,t)$
and the predicted solution $\hat{u}_{\alpha}(\cdot,t)$ at time $t=0.5$,
for $\alpha=0.3$ (left), $\alpha=0.5$ (center) and $\alpha=1.3$
(right).}
\label{level1_2}
\end{figure}

\begin{figure}[H]
\centering{}\includegraphics[scale=1.35]{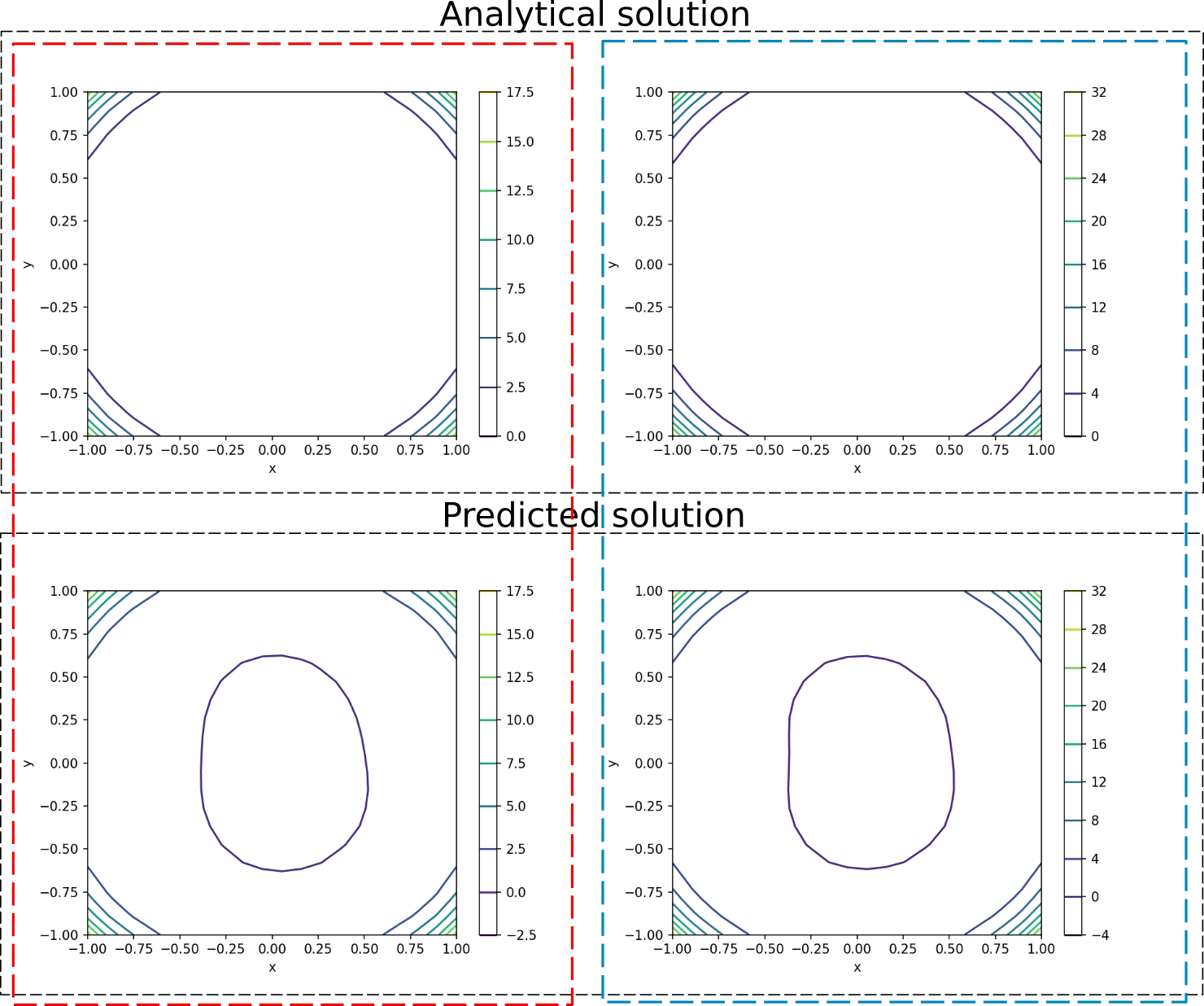}\caption{Level curves of the exact solution $u_{\alpha}(\cdot,t)$ (above)
and the predicted solution $\hat{u}_{\alpha}(\cdot,t)$ (below) for
$\alpha=5$, at times $t=0.5$ (left) and $t=1$ (right).}
\label{level2_2}
\end{figure}

\noindent $\hspace*{1em}$Later, we have examined the contour lines
of $\hat{u}_{\alpha}(\cdot,t)$ for $\alpha\in\{0.3,0.5,1.3,5\}$
and for not very large times $t>0$. For every fixed $\alpha\in\{0.3,0.5,1.3\}$,
the level curves of $\hat{u}_{\alpha}(\cdot,t)$ overlap quite well
those of $u_{\alpha}(\cdot,t$), with some small differences between
one case and the other. More precisely, for each $\alpha\in\{0.5,1.3\}$
the contour lines corresponding to the same level are almost indistinguishable,
at least for not very long times $t$ (see, for example, Figure \ref{level1_2},
where $t=0.5$).\\
For $\alpha=5$ and $t>0$, we have also noted that the approximate
solution is essentially equal to zero in a fairly large region $\varSigma_{t}$
around the origin $(0,0)$ of the $xy$-plane (see Fig. \ref{level2_2}).
As already said for problem (\ref{eq:P1}), this means that such region
consists of numerical zeros of $\hat{u}_{5}(\cdot,t)$, while for
$t>0$ we know that $u_{5}(x,y,t)=0$ if and only if $(x,y)=(0,0)$.
However, this discrepancy is reasonably small for short times, since
the order of magnitude of $u_{5}(x,y,t)$ does not exceed $10^{-2}$
within $\varSigma_{t}$ for $0<t\leq10$. 
 
\noindent $\hspace*{1em}$To evaluate in a more quantitative manner
the accuracy of our method in solving problem (\ref{eq:P2}) and the
numerical convergence of the solution $\hat{u}$ toward $u$, we may
now consider again the quantities (\ref{eq:L2err}) and (\ref{eq:err_rel}).
Passing from Cartesian to polar coordinates, we find 
\[
\|u_{\alpha}(\cdot,t)\|_{L^{2}(\Omega)}^{2}\,=\,t^{2}\underset{\Omega\,\,\,\,}{\iint}(x^{2}+y^{2})^{2\alpha}\,dx\,dy\,=\,\frac{\pi\,t^{2}}{2\alpha+1}\,,
\]
so that we now have
\[
\mathcal{E}_{rel}(T)\,:=\,\,\frac{\mathcal{E}(T)}{\underset{t\,\in\,(0,T)}{\sup}\,\|u_{\alpha}(\cdot,t)\|_{L^{2}(\Omega)}^{2}}\,\,=\,\,\frac{2\alpha+1}{\pi\,T^{2}}\,\underset{t\,\in\,(0,T)}{\sup}\,\|\hat{u}_{\alpha}(\cdot,t)-u_{\alpha}(\cdot,t)\|_{L^{2}(\Omega)}^{2}.
\]
During the experimental phase, we have estimated $\mathcal{E}(T)$ and
$\mathcal{E}_{rel}(T)$ for $\alpha\in\{0.3,0.5,1.3,5\}$ and $T\in\{1,10,20,40,100\}$.
The results that we have obtained are reported in Tables \ref{alpha0.3}$-$\ref{tablealpha5}
and show that, for any fixed value of $\alpha$, the estimate of $\mathcal{E}(T)$ follows an increasing trend, as prescribed by (\ref{eq:L2err}). Furthermore, by analyzing
the orders of magnitude of $\mathcal{E}(T)$ and $\mathcal{E}_{rel}(T)$,
we may affirm that our approach provides very accurate predictions,
on both a short and long-term scale. 
\begin{table}[H]
\noindent \centering{}%
\begin{tabular}{|c|c|c|}
\hline 
\multicolumn{1}{|c}{} & \multicolumn{1}{c}{$\boldsymbol{\alpha=0.3}$} & \tabularnewline
\hline 
\hline 
\textbf{Final time $\boldsymbol{T}$} & \textbf{Estimate of} $\boldsymbol{\mathcal{E}(T)}$ & \textbf{Estimate of} $\boldsymbol{\mathcal{E}_{rel}(T)}$\tabularnewline
\hline 
$1$ & $2.42\times10^{-5}$ & $1.23\times10^{-5}$\tabularnewline
\hline 
$10$ & $4.23\times10^{-4}$ & $2.15\times10^{-6}$\tabularnewline
\hline 
$20$ & $2.20\times10^{-3}$ & $2.80\times10^{-6}$\tabularnewline
\hline 
$40$ & $1.48\times10^{-2}$ & $4.71\times10^{-6}$\tabularnewline
\hline 
$100$ & $2.54\times10^{-1}$ & $1.29\times10^{-5}$\tabularnewline
\hline 
\end{tabular}\caption{Estimates of\textcolor{red}{{} }$\mathcal{E}(T)$ and $\mathcal{E}_{rel}(T)$
for $\alpha=0.3$ and $T\in\{1,10,20,40,100\}$.}
\label{alpha0.3}
\end{table}
 
\begin{table}[H]
\noindent \centering{}%
\begin{tabular}{|c|c|c|}
\hline 
\multicolumn{1}{|c}{} & \multicolumn{1}{c}{$\boldsymbol{\alpha=0.5}$} & \tabularnewline
\hline 
\hline 
\textbf{Final time $\boldsymbol{T}$} & \textbf{Estimate of} $\boldsymbol{\mathcal{E}(T)}$ & \textbf{Estimate of} $\boldsymbol{\mathcal{E}_{rel}(T)}$\tabularnewline
\hline 
$1$ & $2.81\times10^{-5}$ & $1.79\times10^{-5}$\tabularnewline
\hline 
$10$ & $3.41\times10^{-4}$ & $2.17\times10^{-6}$\tabularnewline
\hline 
$20$ & $1.26\times10^{-3}$ & $2.01\times10^{-6}$\tabularnewline
\hline 
$40$ & $9.41\times10^{-2}$ & $3.74\times10^{-5}$\tabularnewline
\hline 
$100$ & $2.25\times10^{-1}$ & $1.43\times10^{-5}$\tabularnewline
\hline 
\end{tabular}\caption{Estimates of\textcolor{red}{{} }$\mathcal{E}(T)$ and $\mathcal{E}_{rel}(T)$
for $\alpha=0.5$ and $T\in\{1,10,20,40,100\}$.}
\label{alpha0.5}
\end{table}
 
\begin{table}[H]
\noindent \centering{}%
\begin{tabular}{|c|c|c|}
\hline 
\multicolumn{1}{|c}{} & \multicolumn{1}{c}{$\boldsymbol{\alpha=1.3}$} & \tabularnewline
\hline 
\hline 
\textbf{Final time $\boldsymbol{T}$} & \textbf{Estimate of} $\boldsymbol{\mathcal{E}(T)}$ & \textbf{Estimate of} $\boldsymbol{\mathcal{E}_{rel}(T)}$\tabularnewline
\hline 
$1$ & $2.85\times10^{-5}$ & $3.27\times10^{-5}$\tabularnewline
\hline 
$10$ & $1.76\times10^{-4}$ & $2.02\times10^{-6}$\tabularnewline
\hline 
$20$ & $1.26\times10^{-3}$ & $3.61\times10^{-6}$\tabularnewline
\hline 
$40$ & $1.51\times10^{-3}$ & $1.08\times10^{-6}$\tabularnewline
\hline 
$100$ & $1.36\times10^{-2}$ & $1.56\times10^{-6}$\tabularnewline
\hline 
\end{tabular}\caption{Estimates of\textcolor{red}{{} }$\mathcal{E}(T)$ and $\mathcal{E}_{rel}(T)$
for $\alpha=1.3$ and $T\in\{1,10,20,40,100\}$.}
\label{alpha1.3}
\end{table}
 
\begin{table}[H]
\noindent \centering{}%
\begin{tabular}{|c|c|c|}
\hline 
\multicolumn{1}{|c}{} & \multicolumn{1}{c}{$\boldsymbol{\alpha=5}$} & \tabularnewline
\hline 
\hline 
\textbf{Final time $\boldsymbol{T}$} & \textbf{Estimate of} $\boldsymbol{\mathcal{E}(T)}$ & \textbf{Estimate of} $\boldsymbol{\mathcal{E}_{rel}(T)}$\tabularnewline
\hline 
$1$ & $4.09\times10^{-4}$ & $1.43\times10^{-3}$\tabularnewline
\hline 
$10$ & $1.01\times10^{-3}$ & $3.53\times10^{-5}$\tabularnewline
\hline 
$20$ & $1.19\times10^{-1}$ & $1.04\times10^{-3}$\tabularnewline
\hline 
$40$ & $7.51\times10^{-1}$ & $1.64\times10^{-3}$\tabularnewline
\hline 
$100$ & $8.78\times10^{-1}$ & $3.07\times10^{-4}$\tabularnewline
\hline 
\end{tabular}\caption{Estimates of\textcolor{red}{{} }$\mathcal{E}(T)$ and $\mathcal{E}_{rel}(T)$
for $\alpha=5$ and $T\in\{1,10,20,40,100\}$.}
\label{tablealpha5}
\end{table}

\subsection{Third test problem}

\begin{figure}[H]
\centering{}\includegraphics[scale=0.85]{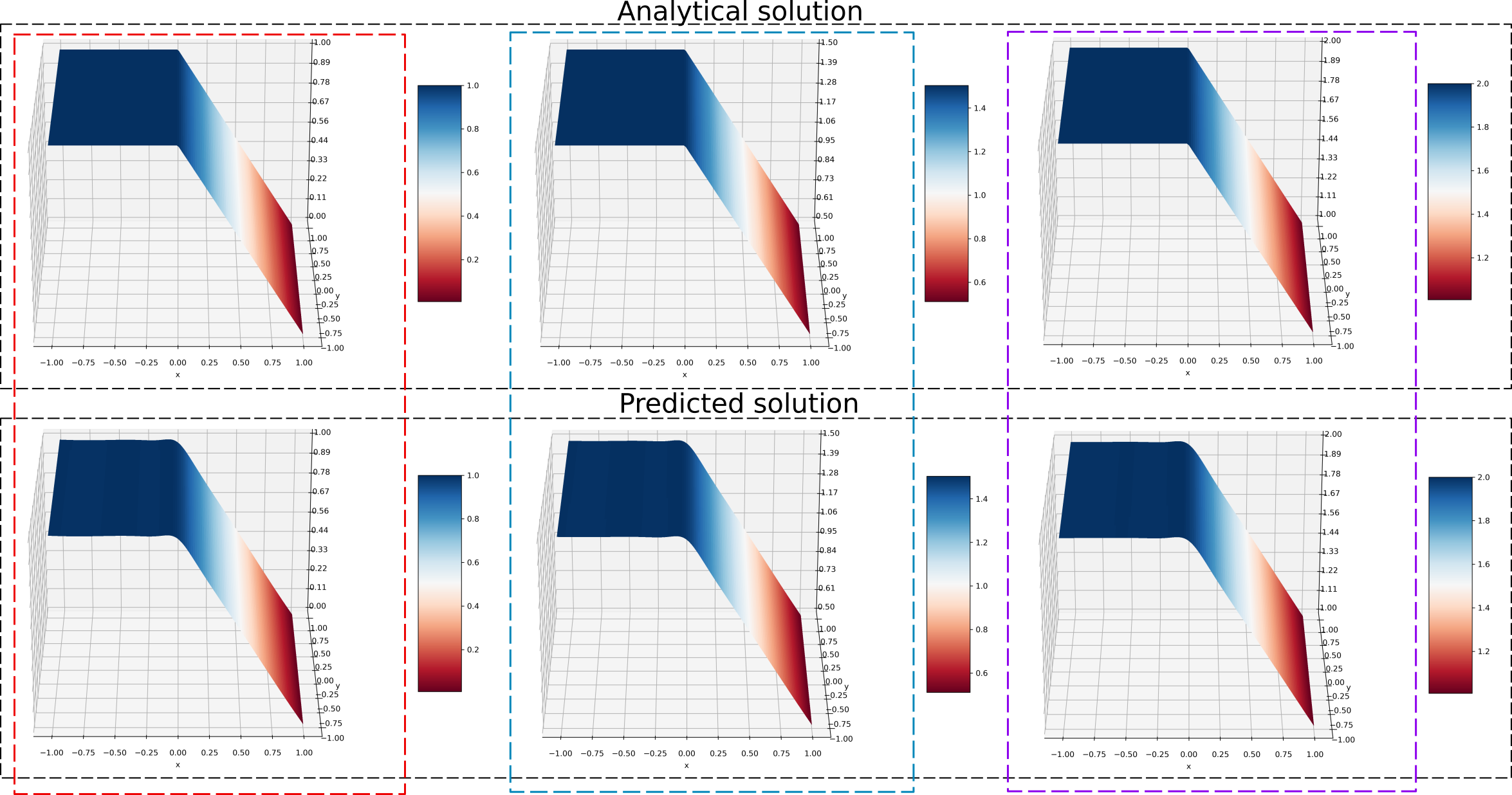}\caption{Plots of the exact solution to problem (\ref{eq:P3}) (above) and
the predicted solution $\hat{u}(\cdot,t)$ (below) for $\textcolor{red}{t=0}$
(left), $\textcolor[RGB]{0,134,185}{t=0.5}$ (center) and $\textcolor[RGB]{141,0,237}{t=1}$
(right).}
\label{pro3_1}
\end{figure}

\noindent $\hspace*{1em}$We shall now consider the problem
\begin{equation}
\begin{cases}
\begin{array}{cc}
\partial_{t}v-\mathrm{div}\left((\vert\nabla v\vert-1)_{+}\,\frac{\nabla v}{\vert\nabla v\vert}\right)=1 & \mathrm{in\,\,}\Omega_{T},\,\,\,\,\,\,\,\,\,\,\,\,\,\,\,\,\,\,\,\,\,\,\,\,\,\,\,\,\,\,\,\,\,\,\,\,\,\,\,\,\,\,\,\,\,\,\,\,\,\\
v(x,y,0)=g(x,y)\,\,\,\,\, & \mathrm{if}\,\,(x,y)\in\overline{\Omega},\,\,\,\,\,\,\,\,\,\,\,\,\,\,\,\,\,\,\,\,\,\,\,\,\,\,\,\,\,\,\,\,\\
v(x,y,t)=h(x,y,t) & \,\,\,\mathrm{if}\,\,(x,y)\in\partial\Omega\,\,\mathrm{\land}\,\,t\in(0,T),
\end{array} & \tag{P3}\end{cases}\label{eq:P3}
\end{equation}
\\
where $\Omega=(-1,1)\times(-1,1)$, 
\[
{\normalcolor {\color{blue}{\normalcolor g(x,y):=\begin{cases}
\begin{array}{cc}
1 & \,\,\,\mathrm{if}\,\,-1\leq x\leq0\,\land\,-1\leq y\leq1,\\
1-x & \mathrm{if}\,\,0<x\leq1\,\land\,-1\leq y\leq1,\,\,\,\,
\end{array} & {\normalcolor }\end{cases}}}}
\]
and
\[
h(x,y,t):=\begin{cases}
\begin{array}{cc}
1+t & \mathrm{if}\,\,-1\leq x\leq0\,\land\,y=-1,\\
1+t & \mathrm{if}\,\,-1\leq x\leq0\,\land\,y=1,\,\,\,\,\,\\
1+t & \mathrm{if}\,\,x=-1\,\land\,-1\leq y\leq1,\,\,\,\\
t & \mathrm{if}\,\,x=1\,\land\,-1\leq y\leq1,\,\,\,\,\,\,\,\\
1-x+t & \mathrm{if}\,\,0<x\leq1\,\land\,y=-1,\,\,\,\,\,\,\,\\
1-x+t & \mathrm{if}\,\,0<x\leq1\,\land\,y=1.\,\,\,\,\,\,\,\,\,\,\,\,
\end{array}\end{cases}
\]
\\
The exact solution of this problem is given
by
\begin{equation}
u(x,y,t)=\begin{cases}
\begin{array}{cc}
1+t & \,\,\,\mathrm{if}\,\,-1\leq x\leq0\,\land\,-1\leq y\leq1,\\
1-x+t & \mathrm{if}\,\,0<x\leq1\,\land\,-1\leq y\leq1.\,\,\,\,
\end{array} & {\normalcolor }\end{cases}\label{eq:sol3}
\end{equation}
Therefore, for any fixed time $t\geq0$, the graph of the function
$u(\cdot,t)$ is given by the union of the horizontal region
\[
\mathcal{H}_{t}:=\left\{ (x,y,1+t):-1\leq x\leq0,\,-1\leq y\leq1\right\} 
\]
and the sliding plane
\[
\mathcal{I}_{t}:=\left\{ (x,y,1-x+t):\,0\leq x\leq1,\,-1\leq y\leq1\right\} .
\]
Let us denote by $\mathcal{G}_{t}:=\mathcal{H}_{t}\cup\mathcal{I}_{t}$
the graph of $u(\cdot,t)$. Then, as time goes by, the set $\mathcal{G}_{t}$
slides along a vertical axis with a constant velocity and no deformation, since $\partial_{t}u\equiv1$ over $\Omega_{T}$ (see
Fig. \ref{pro3_1}, above). 

\begin{figure}[H]
\centering{}\includegraphics[scale=0.49]{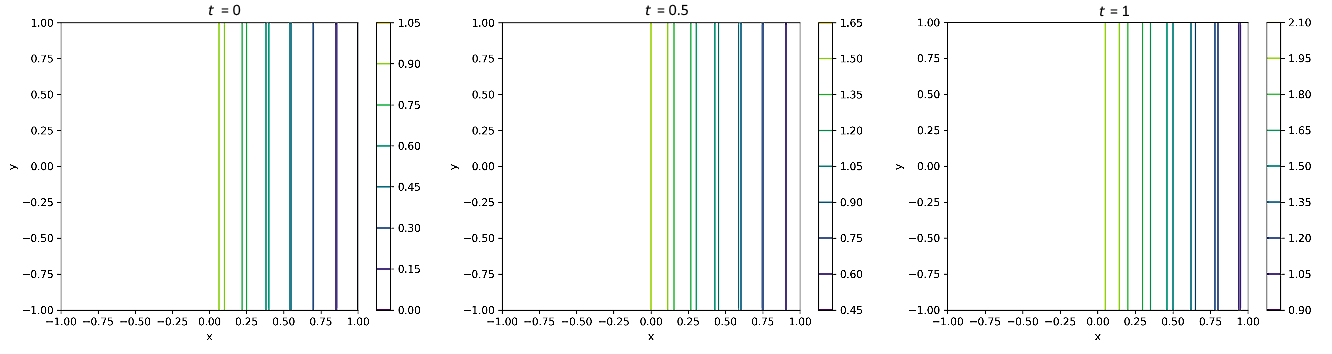}\caption{Superposition of the level curves of the exact solution $u(\cdot,t)$
and the predicted solution $\hat{u}(\cdot,t)$ for $t=0$ (left),
$t=0.5$ (center) and $t=1$ (right).}
\label{prob3_2}
\end{figure}

\noindent $\hspace*{1em}$The plot of the approximate solution $\hat{u}(\cdot,t)$
roughly resembles that of $u(\cdot,t)$ for both short and long times
$t\geq0$, except near the joining line $\mathcal{H}_{t}\cap\mathcal{I}_{t}$,
where the graph of the solution appears to have been slightly smoothed
(see Figure \ref{pro3_1}, below). However, this is not surprising
at all, since we already know that, for any fixed $t\geq0$, the function
$u(\cdot,t):\overline{\Omega}\rightarrow\mathbb{R}$ defined by (\ref{eq:sol3})
is not differentiable at any point of the open segment $S_{0}:=\left\{ (x,y)\in\Omega:x=0\right\} $.
This fact also has repercussions in the comparison between the level
curves of $u(\cdot,t)$ and $\hat{u}(\cdot,t)$, whose superposition
is far from being perfect on approaching the segment $S_{0}$ from
the right, i.e. for $x>0$ (see Fig. \ref{prob3_2}). \\
$\hspace*{1em}$Furthermore, we have also observed that the evolution of $\hat{u}$
over time accurately reflects the evolution of the set $\mathcal{G}_{t}$ described above.

\noindent $\hspace{1em}$In order to assess in a more quantitative way the
accuracy of our method in solving (\ref{eq:P3}) and the distance
between the solutions $u$ and $\hat{u}$, we resort again to the
quantities defined in (\ref{eq:L2err}) and (\ref{eq:err_rel}). Through an easy
calculation, we get 
\[
\|u(\cdot,t)\|_{L^{2}(\Omega)}^{2}\,=\,4\,t^{2}+6\,t+\frac{8}{3}\,,
\]
so that we now have 
\begin{equation}\label{eq:ratio1}
\mathcal{E}_{rel}(T)\,:=\,\,\frac{\mathcal{E}(T)}{\underset{t\,\in\,(0,T)}{\sup}\,\|u(\cdot,t)\|_{L^{2}(\Omega)}^{2}}\,\,=\,\,\frac{3}{12\,T^{2}+18\,T+8}\,\,\,\underset{t\,\in\,(0,T)}{\sup}\,\|\hat{u}(\cdot,t)-u(\cdot,t)\|_{L^{2}(\Omega)}^{2}.
\end{equation}
Proceeding as for the previous problems, we have estimated $\mathcal{E}(T)$
and $\mathcal{E}_{rel}(T)$ for
\[
T\in\{1,10,100,200,300\}.
\]
Table \ref{table_prob3} contains the results obtained and reveals that the estimate of $\mathcal{E}(T)$
exhibits again an increasing behavior, as expected from
(\ref{eq:L2err}). Furthermore, from this table, it seems that the asymptotic trend of the estimate of $\mathcal{E}_{rel}(T)$ may encounter a sort of plateau at $T=100$, after which convergence sensibly slows down. We do not know whether this is a typical behavior, since we cannot draw information from (\ref{eq:ratio1}) in this sense. In fact, from the definition of $\mathcal{E}_{rel}(T)$ it is not possible to predict what the combined effect of $\mathcal{E}(T)$ and $\underset{t\,\in\,(0,T)}{\sup}\,\|u(\cdot,t)\|_{L^{2}(\Omega)}^{2}$ is, since $\mathcal{E}_{rel}(T)$ is the ratio of two functions which are both increasing with respect to $T$ and we cannot determine \textit{a priori} the growth rate of $\mathcal{E}(T)$. Nevertheless, by carefully examining the orders of magnitude
of both $\mathcal{E}(T)$ and $\mathcal{E}_{rel}(T)$, we can conclude
that our method produces accurate results also in this case, in both
short and long-term predictions.
\begin{table}[H]
\noindent \centering{}%
\begin{tabular}{|c|c|c|}
\hline 
\textbf{Final time $\boldsymbol{T}$} & \textbf{Estimate of} $\boldsymbol{\mathcal{E}(T)}$ & \textbf{Estimate of} $\boldsymbol{\mathcal{E}_{rel}(T)}$\tabularnewline
\hline 
$1$ & $1.67\times10^{-1}$ & $1.32\times10^{-2}$\tabularnewline
\hline 
$10$ & $1.69\times10^{-1}$ & $3.65\times10^{-4}$\tabularnewline
\hline 
$100$ & $3.72\times10^{-1}$ & $9.16\times10^{-6}$\tabularnewline
\hline 
$200$ & $6.92\times10^{-1}$ & $4.29\times10^{-6}$\tabularnewline
\hline 
$300$ & $7.21\times10^{-1}$ & $1.99\times10^{-6}$\tabularnewline
\hline 
\end{tabular}\caption{Estimates of\textcolor{red}{{} }$\mathcal{E}(T)$ and $\mathcal{E}_{rel}(T)$
for $T\in\{1,10,100,200,300\}$.}
\label{table_prob3}
\end{table}

\subsection{Fourth test problem}
\noindent $\hspace*{1em}$We now consider the problem
\begin{equation}
\begin{cases}
\begin{array}{cc}
\partial_{t}v-\mathrm{div}\left((\vert\nabla v\vert-1)_{+}\,\frac{\nabla v}{\vert\nabla v\vert}\right)=1 & \mathrm{in\,\,}\Omega_{T},\,\,\,\,\,\,\,\,\,\,\,\,\,\,\,\,\,\,\,\,\,\,\,\,\,\,\,\,\,\,\,\,\,\,\,\,\,\,\,\,\,\,\,\,\,\,\,\,\,\,\,\,\,\\
v(x,y,z,0)=\frac{1}{2}(x^{2}+y^{2}+z^{2}) & \,\,\,\,\,\mathrm{if}\,\,(x,y,z)\in\overline{\Omega},\,\,\,\,\,\,\,\,\,\,\,\,\,\,\,\,\,\,\,\,\,\,\,\,\,\,\,\,\,\,\,\,\,\,\,\\
v(x,y,z,t)=\frac{1}{2}+t\,\,\,\,\,\,\,\,\,\,\,\,\,\,\,\,\,\,\,\,\,\,\,\,\,\, & \,\,\,\,\,\,\mathrm{if}\,\,(x,y,z)\in\partial\Omega\,\,\mathrm{\land}\,\,t\in(0,T),
\end{array} & \tag{P4}\end{cases}\label{eq:P4}
\end{equation}
\\
where $\Omega=\{(x,y,z)\in\mathbb{R}^{3}:x^{2}+y^{2}+z^{2}<1\}$.
This problem is the 3-dimensional version of problem (\ref{eq:P1})
and its exact solution is given by
\[
u(x,y,z,t)=\frac{1}{2}(x^{2}+y^{2}+z^{2})+t.
\]
$\hspace*{1em}$To evaluate the accuracy of our method in solving problem
(\ref{eq:P4}) and the distance between the predicted solution $\hat{u}$
and the exact solution $u$, we confined ourselves to considering
the quantities (\ref{eq:L2err}) and (\ref{eq:err_rel}). Passing
from Cartesian to spherical coordinates, one can easily find that
\[
\|u(\cdot,t)\|_{L^{2}(\Omega)}^{2}\,=\underset{\Omega\,\,\,}{\iiint}\left[\frac{1}{2}(x^{2}+y^{2}+z^{2})+t\right]^{2}dx\,dy\,dz\,=\,\pi\left(\frac{4}{3}\,t^{2}+\frac{4}{5}\,t+\frac{1}{7}\right),
\]
and therefore 
\[
\mathcal{E}_{rel}(T)\,:=\,\,\frac{\mathcal{E}(T)}{\underset{t\,\in\,(0,T)}{\sup}\,\|u(\cdot,t)\|_{L^{2}(\Omega)}^{2}}\,\,=\,\,\frac{105}{\pi\,(140\,T^{2}+84\,T+15)}\,\,\,\underset{t\,\in\,(0,T)}{\sup}\,\|\hat{u}(\cdot,t)-u(\cdot,t)\|_{L^{2}(\Omega)}^{2}.
\]
Proceeding as for problem (\ref{eq:P1}), we have estimated both
$\mathcal{E}(T)$ and $\mathcal{E}_{rel}(T)$ for 
\[
T\in\{1,10,20,30,40,50,100\}.
\]
The data that we have obtained are reported in Table \ref{prob4}
and show that the estimate of $\mathcal{E}(T)$ is monotonically increasing, in agreement with the definition (\ref{eq:L2err}). From Table \ref{prob4} it also emerges that the trend of the estimate of $\mathcal{E}_{rel}(T)$ has a sort of plateau between $T=30$ and $T=40$, after which there is a slight rise. In this regard, the same considerations made for Table \ref{table_prob3} apply. However, for every
$T\leq100$ the order of magnitude of $\mathcal{E}_{rel}(T)$ is not
greater than $10^{-5}$. Therefore, we may affirm that our predictive
method is very accurate and efficient in this case, on both a short
and long time scale.\medskip{}
\begin{table}[H]
\noindent \centering{}%
\begin{tabular}{|c|c|c|}
\hline 
\textbf{Final time $\boldsymbol{T}$} & \textbf{Estimate of} $\boldsymbol{\mathcal{E}(T)}$ & \textbf{Estimate of} $\boldsymbol{\mathcal{E}_{rel}(T)}$\tabularnewline
\hline 
$1$ & $6.59\times10^{-4}$ & $9.21\times10^{-5}$\tabularnewline
\hline 
$10$ & $1.21\times10^{-3}$ & $2.73\times10^{-6}$\tabularnewline
\hline 
$20$ & $1.73\times10^{-3}$ & $1.00\times10^{-6}$\tabularnewline
\hline 
$30$ & $2.90\times10^{-3}$ & $7.54\times10^{-7}$\tabularnewline
\hline 
$40$ & $6.13\times10^{-3}$ & $9.02\times10^{-7}$\tabularnewline
\hline 
$50$ & $4.89\times10^{-2}$ & $4.61\times10^{-6}$\tabularnewline
\hline 
$100$ & $1.44\times10^{-1}$ & $3.42\times10^{-6}$\tabularnewline
\hline 
\end{tabular}\caption{Estimates of $\mathcal{E}(T)$ and $\mathcal{E}_{rel}(T)$
for $T\in\{1,10,20,30,40,50,100\}$.}
\label{prob4}
\end{table}

\subsection{Fifth test problem}
\noindent $\hspace*{1em}$Let $\alpha>0$ and $\omega=2\alpha\,(2\alpha+1)$.
The last problem that we consider is\\
\\
\begin{equation}
\begin{cases}
\begin{array}{cc}
\partial_{t}v-\mathrm{div}\left((\vert\nabla v\vert-1)_{+}\,\frac{\nabla v}{\vert\nabla v\vert}\right)=f & \mathrm{in\,\,}\Omega_{T},\,\,\,\,\,\,\,\,\,\,\,\,\,\,\,\,\,\,\,\,\,\,\,\,\,\,\,\,\,\,\,\,\,\,\,\,\,\,\,\,\,\,\,\,\,\,\,\,\,\,\,\,\,\,\,\,\\
v(x,y,z,0)=0 & \,\,\,\mathrm{if}\,\,(x,y,z)\in\overline{\Omega},\,\,\,\,\,\,\,\,\,\,\,\,\,\,\,\,\,\,\,\,\,\,\,\,\,\,\,\,\,\,\,\,\,\,\,\,\\
v(x,y,z,t)=t & \,\,\,\mathrm{if}\,\,(x,y,z)\in\partial\Omega\,\,\mathrm{\land}\,\,t\in(0,T),
\end{array} & \tag{P5}\end{cases}\label{eq:P5}
\end{equation}
\\
\\
where $\Omega=\{(x,y,z)\in\mathbb{R}^{3}:x^{2}+y^{2}+z^{2}<1\}$
and\\
\\ 
\[
{\normalcolor {\color{blue}{\normalcolor f(x,y,z,t):=\begin{cases}
\begin{array}{cc}
(x^{2}+y^{2}+z^{2})^{\alpha} & \mathrm{if}\,\,2\alpha\,t\,(x^{2}+y^{2}+z^{2})^{\alpha-\frac{1}{2}}\leq1,\\
(x^{2}+y^{2}+z^{2})^{\alpha-1}(x^{2}+y^{2}+z^{2}-\omega\,t)+\,\frac{2}{\sqrt{x^{2}+y^{2}+z^{2}}} & \mathrm{if}\,\,2\alpha\,t\,(x^{2}+y^{2}+z^{2})^{\alpha-\frac{1}{2}}>1.
\end{array}\end{cases}}}}
\]
\\
This problem is nothing but the 3-dimensional
version of (\ref{eq:P2}) and its exact solution is given by
\[
u(x,y,z,t)\equiv u_{\alpha}(x,y,z,t):=t\,(x^{2}+y^{2}+z^{2})^{\alpha}.
\]
$\hspace*{1em}$In order to assess the accuracy of our method in solving (\ref{eq:P5}) and the distance between the approximate solution
$\hat{u}$ and the exact solution $u$, we again limited ourselves to estimating
the quantities defined in (\ref{eq:L2err}) and (\ref{eq:err_rel}). Switching
from Cartesian to spherical coordinates, we can easily obtain 
\[
\|u_{\alpha}(\cdot,t)\|_{L^{2}(\Omega)}^{2}\,=\,t^{2}\underset{\Omega\,\,\,}{\iiint}(x^{2}+y^{2}+z^{2})^{2\alpha}\,dx\,dy\,dz\,=\,\frac{4\,\pi\,t^{2}}{4\,\alpha+3}\,.
\]
This yields
\[
\mathcal{E}_{rel}(T)\,=\,\,\frac{4\,\alpha+3}{4\,\pi\,T^{2}}\,\,\underset{t\,\in\,(0,T)}{\sup}\,\|\hat{u}(\cdot,t)-u(\cdot,t)\|_{L^{2}(\Omega)}^{2}.
\]
Proceeding as for problem (\ref{eq:P2}), we have estimated $\mathcal{E}(T)$
and $\mathcal{E}_{rel}(T)$ for $\alpha\in\{0.3,0.5,1.3,5\}$ and
$T\in\{1,10,20,40,100\}$. The results that have been obtained are
shown in Tables \ref{alpha0.3_pro5}$-$\ref{tablealpha5_pro5} and
reveal that, for any fixed value of $\alpha$, the
estimate of $\mathcal{E}(T)$ is again monotonically increasing, as expected from (\ref{eq:L2err}). Nevertheless, by carefully analyzing
the orders of magnitude of both $\mathcal{E}(T)$ and $\mathcal{E}_{rel}(T)$,
we can deduce that our method provides accurate solutions also in
this case, in both short and long-term predictions. 
\begin{table}[H]
\noindent \centering{}%
\begin{tabular}{|c|c|c|}
\hline 
\multicolumn{1}{|c}{} & \multicolumn{1}{c}{$\boldsymbol{\alpha=0.3}$} & \tabularnewline
\hline 
\hline 
\textbf{Final time $\boldsymbol{T}$} & \textbf{Estimate of} $\boldsymbol{\mathcal{E}(T)}$ & \textbf{Estimate of} $\boldsymbol{\mathcal{E}_{rel}(T)}$\tabularnewline
\hline 
$1$ & $1.27\times10^{-5}$ & $4.25\times10^{-6}$\tabularnewline
\hline 
$10$ & $6.93\times10^{-4}$ & $2.31\times10^{-6}$\tabularnewline
\hline 
$20$ & $4.64\times10^{-3}$ & $3.88\times10^{-6}$\tabularnewline
\hline 
$40$ & $2.37\times10^{-2}$ & $4.95\times10^{-6}$\tabularnewline
\hline 
$100$ & $2.58\times10^{-1}$ & $8.63\times10^{-6}$\tabularnewline
\hline 
\end{tabular}\caption{Estimates of $\mathcal{E}(T)$ and $\mathcal{E}_{rel}(T)$
for $\alpha=0.3$ and $T\in\{1,10,20,40,100\}$.}
\label{alpha0.3_pro5}
\end{table}
 
\begin{table}[H]
\noindent \centering{}%
\begin{tabular}{|c|c|c|}
\hline 
\multicolumn{1}{|c}{} & \multicolumn{1}{c}{$\boldsymbol{\alpha=0.5}$} & \tabularnewline
\hline 
\hline 
\textbf{Final time $\boldsymbol{T}$} & \textbf{Estimate of} $\boldsymbol{\mathcal{E}(T)}$ & \textbf{Estimate of} $\boldsymbol{\mathcal{E}_{rel}(T)}$\tabularnewline
\hline 
$1$ & $1.02\times10^{-5}$ & $4.05\times10^{-6}$\tabularnewline
\hline 
$10$ & $4.60\times10^{-4}$ & $1.83\times10^{-6}$\tabularnewline
\hline 
$20$ & $2.25\times10^{-3}$ & $2.24\times10^{-6}$\tabularnewline
\hline 
$40$ & $3.44\times10^{-2}$ & $8.54\times10^{-6}$\tabularnewline
\hline 
$100$ & $4.47\times10^{-1}$ & $1.78\times10^{-5}$\tabularnewline
\hline 
\end{tabular}\caption{Estimates of $\mathcal{E}(T)$ and $\mathcal{E}_{rel}(T)$
for $\alpha=0.5$ and $T\in\{1,10,20,40,100\}$.}
\label{alpha0.5_pro5}
\end{table}
 
\begin{table}[H]
\noindent \centering{}%
\begin{tabular}{|c|c|c|}
\hline 
\multicolumn{1}{|c}{} & \multicolumn{1}{c}{$\boldsymbol{\alpha=1.3}$} & \tabularnewline
\hline 
\hline 
\textbf{Final time $\boldsymbol{T}$} & \textbf{Estimate of} $\boldsymbol{\mathcal{E}(T)}$ & \textbf{Estimate of} $\boldsymbol{\mathcal{E}_{rel}(T)}$\tabularnewline
\hline 
$1$ & $3.21\times10^{-5}$ & $2.09\times10^{-5}$\tabularnewline
\hline 
$10$ & $9.83\times10^{-4}$ & $6.42\times10^{-6}$\tabularnewline
\hline 
$20$ & $9.73\times10^{-3}$ & $1.59\times10^{-5}$\tabularnewline
\hline 
$40$ & $2.61\times10^{-2}$ & $1.06\times10^{-5}$\tabularnewline
\hline 
$100$ & $1.56\times10^{-1}$ & $1.02\times10^{-5}$\tabularnewline
\hline 
\end{tabular}\caption{Estimates of $\mathcal{E}(T)$ and $\mathcal{E}_{rel}(T)$
for $\alpha=1.3$ and $T\in\{1,10,20,40,100\}$.}
\label{alpha1.3_pro5}
\end{table}
 
\begin{table}[H]
\noindent \centering{}%
\begin{tabular}{|c|c|c|}
\hline 
\multicolumn{1}{|c}{} & \multicolumn{1}{c}{$\boldsymbol{\alpha=5}$} & \tabularnewline
\hline 
\hline 
\textbf{Final time $\boldsymbol{T}$} & \textbf{Estimate of} $\boldsymbol{\mathcal{E}(T)}$ & \textbf{Estimate of} $\boldsymbol{\mathcal{E}_{rel}(T)}$\tabularnewline
\hline 
$1$ & $9.72\times10^{-4}$ & $1.78\times10^{-3}$\tabularnewline
\hline 
$10$ & $1.19\times10^{-3}$ & $2.19\times10^{-5}$\tabularnewline
\hline 
$20$ & $3.61\times10^{-2}$ & $1.65\times10^{-4}$\tabularnewline
\hline 
$40$ & $2.45\times10^{-1}$ & $2.80\times10^{-4}$\tabularnewline
\hline 
$100$ & $6.83\times10^{-1}$ & $1.25\times10^{-4}$\tabularnewline
\hline 
\end{tabular}\caption{Estimates of $\mathcal{E}(T)$ and $\mathcal{E}_{rel}(T)$
for $\alpha=5$ and $T\in\{1,10,20,40,100\}$.}
\label{tablealpha5_pro5}
\end{table}

\section{Conclusions}\label{sec:conclusions}
\noindent $\hspace*{1em}$In this paper, we have explored the ability
of PINNs to accurately predict the solutions of some strongly degenerate
parabolic problems arising in gas filtration through porous media.
Since there are no general methods for finding analytical solutions
to such problems, it is essential to use efficient and accurate numerical
methods. blueOne of the most prevalent methods for addressing these problems is the Finite Difference Method (FDM), wherein PDEs are discretized into a system of algebraic equations to be solved numerically. However, the FDM necessitates the discretization of the domain into a grid of cells or nodes, which can become computationally expensive for large and intricate systems. Although the primary objective of this article is not to prove the effectiveness of a PINN compared to a classical numerical method, we engaged in a comparison with the FDM. As established in the literature, for problems characterized by a less complex domain, the FDM typically exhibits a higher level of accuracy compared to PINNs. Nevertheless, in our study, the advantage of using a PINN lies in the ability to test the model on various presented problems (varying the initial/boundary functions and the $\alpha$ parameter), once it has been trained. Additionally, the FDM can be utilized as a benchmark in cases where the solution to the problem is unknown, ensuring a fair comparison under equivalent accuracy conditions.\\
$\hspace*{1em}$For the test problems discussed here, whose exact
solutions are fortunately known, we have compared the plots of the
predicted solutions with those of the analytical solutions. Moreover,
to evaluate the accuracy of our predictive method in a purely quantitative
way, we have also analyzed the error trends over time. The proposed
approach provides accurate results in line with expectations, at least
in short-term predictions. However, some issues remain open, such
as how to obtain fully reliable plots for the predicted solution when
the exact (unknown) one is not differentiable somewhere, and how to
reduce or eliminate some slight discrepancies between the contour
lines of the predicted solution and those of the analytical solution
in the case $n=2$.\\
$\hspace*{1em}$To the best of our knowledge, this is one of the first
papers demonstrating the effectiveness of the PINN framework for solving
strongly degenerate parabolic problems with asymptotic structure of
Laplacian type.% COMMENTO PER PIA: Classical numerical methods like finite differences, finite elements and finite volumes have some limitations, including high computational costs as the problem dimensions increase.

\noindent \bigskip{}

\noindent \textbf{Acknowledgements.} We would like to thank the reviewers for their valuable comments, which helped to improve our paper. S. Cuomo also acknowledges GNCS-INdAM and the UMI-TAA, UMI-AI research groups. This work has been supported by the following projects:
\begin{itemize}
    \item The Programma Operativo Nazionale Ricerca e Innovazione $2014$-$2020$ (CCI2014IT16M2OP005) Dottorati e contratti di ricerche su tematiche dell'innovazione XXXVII Ciclo, code DOT1318347, CUP: E65F21003980003.
    \item PNRR Centro Nazionale HPC, Big Data e Quantum Computing, (CN\_00000013)(CUP: E63C22000980007), under the NRRP MUR program funded by the NextGenerationEU.
    \item P. Ambrosio has been partially supported by the INdAM$-$GNAMPA 2023 Project “Risultati di regolarità per PDEs in spazi di funzione non-standard” (CUP: E53C22001930001) and by the INdAM$-$GNAMPA 2024 Project “Fenomeno di Lavrentiev, Bounded Slope Condition e regolarità per minimi di funzionali integrali con crescite non standard e lagrangiane non uniformemente convesse” (CUP: E53C23001670001).
\end{itemize}
\bigskip{}

\noindent \textbf{Declarations.} The authors declare that they have no conflict of interest.

\lyxaddress{\noindent \textbf{$\quad$}\\
$\hspace*{1em}$\textbf{Pasquale Ambrosio}\\
 Dipartimento di Matematica e Applicazioni ``R. Caccioppoli''\\
 Università degli Studi di Napoli ``Federico II''\\
 Via Cintia, 80126 Napoli, Italy.\\
 \textit{E-mail address}: pasquale.ambrosio2@unina.it}

\lyxaddress{\noindent $\hspace*{1em}$\textbf{Salvatore Cuomo}\\
 Dipartimento di Matematica e Applicazioni ``R. Caccioppoli''\\
 Università degli Studi di Napoli ``Federico II''\\
 Via Cintia, 80126 Napoli, Italy.\\
 \textit{E-mail address}: salvatore.cuomo@unina.it}

\lyxaddress{\noindent $\hspace*{1em}$\textbf{Mariapia De Rosa}\\
 Dipartimento di Matematica e Applicazioni ``R. Caccioppoli''\\
 Università degli Studi di Napoli ``Federico II''\\
 Via Cintia, 80126 Napoli, Italy.\\
 \textit{E-mail address}: mariapia.derosa@unina.it}

\end{document}